\documentclass[a4paper,12pt]{article}
\usepackage{graphicx}
\usepackage{float}
\usepackage{amsmath}
\usepackage{amssymb}
\usepackage{amsfonts}
\usepackage{theorem}
\usepackage{amscd}
\usepackage{latexsym}
\usepackage{dsfont}
\usepackage[latin1]{inputenc}

\usepackage{eurosym}

\numberwithin{equation}{section}

\newtheorem{rem}{Remark}[section]
\newtheorem{assum}{Assumption}[section]
\newtheorem{theorem}{Theorem}[section]
\newtheorem{lemma}{Lemma}[section]

\newtheorem{prop}{Proposition}[section]

\newtheorem{definition}{Definition}[section]

\DeclareMathOperator*{\argmin}{argmin}
\DeclareMathOperator*{\argminloc}{argminloc}
\DeclareMathOperator*{\diam}{diam}
\DeclareMathOperator*{\dist}{dist}
\DeclareMathOperator*{\supp}{supp}
\newcommand{\Ind}{\mathds{1}}

\newcommand{\zbf}{{\bf z}}

\newcommand{\real}{\mathbb{R}}
\newcommand{\esp}{\mathbb{E}}
\newcommand{\prob}{\mathbb{P}}

\newcommand{\vertbar}{\left\vert\vphantom{\frac{1}{1}}\right.}
\newcommand{\proofend}{\begin{flushright} $\square$ \end{flushright}}
\newcommand{\wch}{\check{w}}
\newcommand{\Gcal}{\mathcal{G}}
\newcommand{\Dcal}{\mathcal{D}}
\newcommand{\rdkap}{\left(\real^d\right)^\kappa}

\begin{document}
\begin{center}

{\sc \Large Convergence of distributed asynchronous learning vector quantization algorithms. \\
\vspace{0.7cm}}

Beno\^{i}t PATRA $^{\mbox{\footnotesize a,b}}$
\vspace{0.5cm}

$^{\mbox{\footnotesize a}}$ LSTA\\Universit\'e Pierre et Marie Curie -- Paris VI\\
Tour 15-25, 4 place Jussieu\\
75252 Paris cedex 05, France\\
\smallskip
\textsf{benoit.patra@upmc.fr}\\
\bigskip
$^{\mbox{\footnotesize b}}$
LOKAD SAS\\
10 rue Philippe de Champaigne\\
75013 Paris, France\\
\smallskip
\textsf{benoit.patra@lokad.com}\\

\vspace{0.5cm}

\end{center}
\begin{abstract}
\noindent {\rm Motivated by the problem of effectively executing clustering algorithms on very large data sets, we address a model for large scale distributed clustering methods. To this end, we briefly recall some standards on the quantization problem and some results on the almost sure convergence of the Competitive Learning Vector Quantization (CLVQ) procedure. A general model for linear distributed asynchronous algorithms well adapted to several parallel computing architectures is also discussed. Our approach brings together this scalable model and the CLVQ algorithm, and we call the resulting technique the Distributed Asynchronous Learning Vector Quantization algorithm (DALVQ). An in-depth analysis of the almost sure convergence of the DALVQ algorithm is performed. A striking result is that we prove that the multiple versions of the quantizers distributed among the processors in the parallel architecture asymptotically reach a consensus almost surely. Furthermore, we also show that these versions converge almost surely towards the same nearly optimal value for the quantization criterion.
}\\

\medskip
\end{abstract}

\section{Introduction}
Distributed algorithms arise in a wide range of applications, including telecommunications, distributed information processing, scientific computing, real time process control and many others. Parallelization is one of the most promising ways to harness greater computing resources, whereas building faster serial computers is increasingly expensive and also faces some physical limits such as transmission speeds and miniaturization. One of the challenges proposed for Machine Learning is to build scalable applications that quickly process large amounts of data in sophisticated ways. Building such large scale algorithms attacks several problems in a distributed framework, such as communication delays in the network or numerous problems caused by the lack of shared memory.\\

Clustering algorithms are one of the primary tools of unsupervised learning. From a practical perspective, clustering plays an outstanding role in data mining applications such as text mining, web analysis, marketing, medical diagnostics, computational biology and many others. Clustering is a separation of data into groups of similar objects. As clustering represents the data with fewer clusters, there is a necessary loss of certain fine details, but simplification is achieved. The popular Competitive Learning Vector Quantization (CLVQ) algorithm (see Gersho and Gray \cite{GER1}) provides a technique for building reliable clusters characterized by their prototypes. As pointed out by Bottou in \cite{BOT4}, the CLVQ algorithm can also be viewed as the on-line version of the widespread Lloyd's method (see Lloyd's \cite{LL01} for the definition) which is referred to as batch $k$-means in \cite{BOT4}. The CLVQ also belongs to the class of stochastic gradient descent algorithms (for more information on stochastic gradient descent procedures we refer the reader to Benveniste et al.~\cite{BEN1}).\\

The analysis of parallel stochastic gradient procedures in a Machine Learning context has recently received a great deal of attention (see for instance Zinkevich et al.~\cite{ZIN} and Mac Donald et al.~\cite{MacDo}). In the present paper, we go further by introducing a model that brings together the original CLVQ algorithm and the comprehensive theory of asynchronous parallel linear algorithms developed by  Tsitsiklis \cite{TSI2}, Tsitsiklis et al.~\cite{TSI1} and Bertsekas and Tsitsiklis \cite{BER1}. The resulting model will be called Distributed Asynchronous Learning Vector Quantization (DALVQ for short). At a high level, the DALVQ algorithm parallelizes several executions of the CLVQ method concurrently on different processors while the results of these algorithms are broadcast through the distributed framework asynchronously and efficiently. Here, the term processor refers to any computing instance in a distributed architecture (see Bullo et al.~\cite[Chapter 1]{BUL1} for more details). Let us remark that there is a series of publications similar in spirit to this paper. Indeed in Frasca et al.~\cite{FRA1} and in Durham et al.~\cite{DUR3}, a coverage control problem is formulated as an optimization problem where the functional cost to be minimized is the same of the quantization problem stated in this manuscript.\\

Let us provide a brief mathematical introduction to the CLVQ technique and DALVQ algorithms. The first technique computes quantization scheme for $d$ dimensional samples $\zbf_1, \zbf_2, \ldots$ using the following iterations on a $\rdkap$ vector,
\begin{equation*}\label{eq:finiteSample}
w(t+1)= w(t) - \varepsilon_{t+1}H\left(\zbf_{t+1},w(t) \right), \quad t \geq 0.
\end{equation*}
In the equation above, $w(0) \in \left(\real^d\right)^\kappa$ and the $\varepsilon_t$ are positive reals. The vector $H(\zbf,w)$ is the opposite of the difference between the sample $\zbf$ and its nearest component in $w$. Assume that there are $M$ computing entities, the data are split among the memory of these machines: $\zbf^i_1, \zbf^i_2, \ldots$, where $i\in \{1,\ldots,M\}$. Therefore, the DALVQ algorithms are defined by the $M$ iterations $\{w^i(t)\}_{t=0}^\infty$, called versions, satisfying (with slight simplifications)
\begin{equation}\label{eq:defwi}
w^i(t+1) =\sum_{j=1}^M a^{i,j}(t)w^j(\tau^{i,j}(t)) - \varepsilon^i_{t+1} H\left(\zbf^i_{t+1},w^i(t) \right),
\end{equation}
$i \in \left\{1, \ldots, M \right\}$ and $t \geq 0$. The time instants $\tau^{i,j}(t) \geq 0$ are deterministic but unknown and the delays satisfy $t - \tau^{i,j}(t) \geq 0$. The families $\{a^{i,j}(t)\}_{j=1}^M$ define the weights of convex combinations.\\

As a striking result, we prove that multiple versions of the quantizers, distributed among the processors in a parallel architecture, asymptotically reach a consensus almost surely. Using the materials introduced above, it writes
\[w^i(t) -w^j(t) \xrightarrow[t \rightarrow \infty]{} 0, \quad \mbox{$(i,j) \in \left\{ 1, \ldots, M\right\}^2$, almost surely (a.s.).}\]
Furthermore, we also show that these versions converge almost surely towards (the same) nearly optimal value for the quantization criterion. These convergence results are similar in spirit to the most satisfactory almost sure convergence theorem for the CLVQ algorithm obtained by Pagès in \cite{PAG1}.\\

For a given time span, our parallel DALVQ algorithm is able to process much more data than a single processor execution of the CLVQ procedure. Moreover, DALVQ is also asynchronous. This means that local algorithms do not have to wait at preset points for messages to become available. This allows some processors to compute faster and execute more iterations than others, and it also allows communication delays to be substantial and unpredictable. The communication channels are also allowed to deliver messages out of order, that is, in a different order than the one in which they were transmitted. Asynchronism can provide two major advantages. First, a reduction of the synchronization penalty, which could bring a speed advantage over a synchronous execution. Second, for potential industrialization, asynchronism has greater implementation flexibility. Tolerance to system failures and uncertainty can also be increased. As in the case with any on-line algorithm, DALVQ also deals with variable data loads over time. In fact, on-line algorithms avoid tremendous and non scalable batch requests on all data sets. Moreover, with an on-line algorithm, new data may enter the system and be taken into account while the algorithm is already running.\\

The paper is organized as follows. In Section \ref{sect:CLVQ} we review some standard facts on the clustering problem. We extract the relevant material from Pagès \cite{PAG1} without proof, thus making our exposition self-contained. In Section \ref{sect:MDA} we give a brief exposition of the mathematical framework for parallel asynchronous gradient methods introduced by Tsitsiklis et al. in \cite{TSI1} and Bertsekas and Tsitsiklis \cite{TSI1,BER1}. The results of Blondel et al.~\cite{BLO1} on the asymptotic consensus in asynchronous parallel averaging problems are also recalled. In Section \ref{sect:DALVQ}, our main results are stated and proved.

\section{Quantization and CLVQ algorithm}\label{sect:CLVQ}
\subsection{Overview}
Let $\mu$ be a probability measure on $\real^d$ with finite second-order moment. The quantization problem consists in finding a ``good approximation'' of $\mu$ by a set of $\kappa$ vectors of $\real^d$ called quantizer. Throughout the document the $\kappa$ quantization points (or prototypes) will be seen as the components of a $\rdkap$-dimensional vector $w=(w_1, \ldots,w_\kappa)$. To measure the correctness of a quantization scheme given by $w$, one introduces a cost function called distortion, defined by
\[C_\mu(w) = \frac12 \int_{\real^d}{\min_{1 \leq \ell \leq \kappa} \left\|\zbf-w_\ell \right\|^2 d\mu(\zbf)}.\]
Under some minimal assumptions, the existence of an optimal quantizer vector $w^\circ \in \argmin_{w \in \rdkap}{C_\mu(w)}$ has been established by Pollard in \cite{POL2} (see also Sabin and Gray \cite[Appendix 2]{SAB1}).\\

In a statistical context, the distribution $\mu$ is only known through $n$ independent random observations $\zbf_1, \ldots, \zbf_n$ drawn according to $\mu$. Denote by $\mu_n$ the empirical distribution based on $\zbf_1, \ldots, \zbf_n$, that is, for every Borel subset $A$ of $\real^d$
\[\mu_n(A) = \frac1n\sum_{i=1}^n{\mathds{1}_{\left\{\zbf_i \in A\right\}}}.\]
Much attention has been devoted to the convergence study of the quantization scheme provided by the empirical minimizers
\[w_n^\circ \in \argmin_{w \in \rdkap}{C_{\mu_n}(w)}.\]

The almost sure convergence of $C_\mu\left(w_n^\circ\right)$ towards $\min_{w \in \rdkap} C_\mu(w)$ was proved by Pollard in \cite{POL3,POL2} and Abaya and Wise in \cite{ABA1}. Rates of convergence and nonasymptotic performance bounds have been considered by Pollard \cite{POL4}, Chou \cite{CHO1}, Linder et al.~\cite{LIN3}, Bartlett et al.~\cite{BAR1}, Linder \cite{LIN2,LIN4}, Antos \cite{ANT1} and Antos et al.~\cite{ANT2}. Convergence results have been established by Biau et al. in \cite{BIA3} where $\mu$ is a measure on a Hilbert space. It turns out that the minimization of the empirical distortion is a computationally hard problem. As shown by Inaba et al. in \cite{INA1}, the computational complexity of this minimization problem is exponential in the number of quantizers $\kappa$ and the dimension of the data $d$. Therefore, exact computations are untractable for most of the practical applications.\\

Based on this, our goal in this document is to investigate effective methods that produce accurate quantizations with data samples. One of the most popular procedure is Lloyd's algorithm (see Lloyd \cite{LL01}) sometimes refereed to as batch $k$-means. A convergence theorem for this algorithm is provided by Sabin and Gray in \cite{SAB1}. Another celebrated quantization algorithm is the Competitive Learning Vector Quantization (CLVQ), also called on-line $k$-means. The latter acronym outlines the fact that data arrive over time while the execution of the algorithm and their characteristics are unknown until their arrival times. The main difference between the CLVQ and the Lloyd's algorithm is that the latter run in batch training mode. This means that the whole training set is presented before performing an update, whereas the CLVQ algorithm uses each item of the training sequence at each update.\\

The CLVQ procedure can be seen as a stochastic gradient descent algorithm. In the more general context of gradient descent methods, one cannot hope for the convergence of the procedure towards global minimizers with a non convex objective function (see for instance Benveniste et al.~\cite{BEN1}). In our quantization context, the distortion mapping $C_\mu$ is not convex (see for instance Graf and Luschgy \cite{GRA1}). Thus, just as in Lloyd's method, the iterations provided by the CLVQ algorithm converge towards local minima of $C_{\mu}$.\\


Assuming that the distribution $\mu$ has a compact support and a bounded density with respect to the Lebesgue measure, Pag\`es states in \cite{PAG1} a result regarding the almost sure consistency of the CLVQ algorithm towards critical points of the distortion $C_{\mu}$. The author shows that the set of critical points necessarily contains the global and local optimal quantizers. The main difficulties in the proof arise from the fact that the gradient of the distortion is singular on $\kappa$-tuples having equal components and the distortion function $C_{\mu}$ is not convex. This explains why standard theories for stochastic gradient algorithm do not apply in this context.\\

\subsection{The quantization problem, basic properties}\label{sec:VQ:Notations}
In the sequel, we denote by $\Gcal$ the closed convex hull of $\supp\left(\mu \right)$, where  $\supp\left(\mu \right)$ stands for the support of the distribution. Observe that, with this notation, the distortion mapping is the function $C : \left(\real^d \right)^\kappa  \longrightarrow  [0, \infty)$ defined by
\begin{equation}\label{def:C}
 C(w) \triangleq \frac12
 \displaystyle\int_{\Gcal}{\min_{1 \leq \ell \leq \kappa} \left\|\zbf-w_\ell \right\|^2 d\mu(\zbf)}, \quad w=\left(w_1, \ldots,w_\kappa\right) \in \left(\real^d \right)^\kappa.
\end{equation}
Throughout the document, with a slight abuse of notation, $\left\|. \right\|$ means both the Euclidean norm of $\real^d$ or $\rdkap$. In addition, the notation $\Dcal_*^\kappa$ stands for the set of all vector of $\rdkap$ with pairwise distinct components, that is,

\begin{equation*}
\Dcal_*^\kappa \triangleq \left\{w \in \rdkap \; \vert \; w_\ell \neq w_k \text{ if and only if } \ell \neq k \right\}.
\end{equation*}

%

\begin{figure}[h]
\begin{center}
\includegraphics[scale = 0.35]{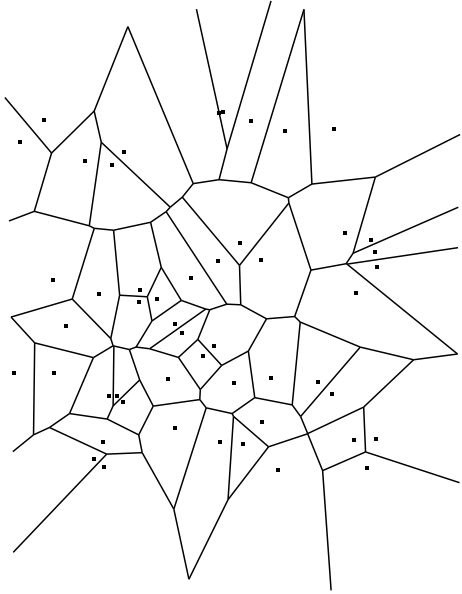}
\end{center}
\caption{Vorono\"i tessellation of 50 points of $\real^2$ drawn uniformly in a square.}
\label{fig:voronoi}
\end{figure}

Under some extra assumptions on $\mu$, the distortion function can be rewritten using space partition set called Vorono\"i tessellation.
\begin{definition}\label{def:tessels}
Let $w \in \rdkap$, the Vorono\"i tessellation of $\Gcal$ related to $w$ is the family of open sets $\left\{W_\ell(w)\right\}_{1 \leq \ell \leq \kappa}$ defined as follows:
\begin{itemize}
\item If $w \in \Dcal_*^\kappa$, for all $1 \leq \ell \leq \kappa$,
\[W_\ell(w) = \left\{v \in \Gcal \; \vertbar \; \left\|w_\ell - v\right\| < \min_{k \neq \ell} \left\|w_k - v\right\|\right\}.\]
\item If $w \in \rdkap \setminus \Dcal_*^\kappa$, for all $1 \leq \ell \leq \kappa$,\\
\begin{itemize}
\item if $\ell = \min\left\{k \; \vert \; w_k=w_\ell\right\}$,
\[W_\ell(w) =\left\{v \in \Gcal \; \vertbar \; \left\|w_\ell - v\right\| < \min_{w_k \neq w_\ell}\left\|w_k - v\right\| \right\} \]
\item otherwise, $W_\ell(w) = \emptyset$.
\end{itemize}
\end{itemize}
\end{definition}
As an illustration, Figure \ref{fig:voronoi} shows Vorono\"i tessellations associated to a vector $w \in ([0,1]\times[0,1])^{50}$ whose components have been drawn independently and uniformly. This figure also highlights a remarkable property of the cell borders, which are portions of hyperplanes (see Graf and Luschgy \cite{GRA1}).\\

Observe that if $\mu(H)$ is zero for any hyperplane $H$ of $\real^d$ (a property which is sometimes referred to as strong continuity) then using Definition \ref{def:tessels}, it is easy to see that the distortion takes the form:
\begin{equation}
C(w) = \frac{1}{2} \sum_{\ell = 1}^\kappa{\displaystyle \int_{W_\ell(w)}{\left\|\zbf - w_\ell \right\|^2 d\mu(\zbf)}}, \quad w \in \rdkap.
\end{equation}

The following assumption will be needed throughout the paper. This assumption is similar to the peak power constraint (see Chou \cite{CHO1} and Linder \cite{LIN4}). Note that most of the results of this subsection are still valid if $\mu$ satisfies the weaker strong continuity property.
\begin{assum}[Compact Supported Density]\label{assum:compactdensity}
The probability measure $\mu$ has a bounded density with respect to the Lebesgue measure on $\real^d$. Moreover, the support of $\mu$ is equal to its convex hull $\Gcal$, which in turn, is compact.
\end{assum}

The next proposition states the differentiability of the distortion $C$, and provides an explicit formula for the gradient $\nabla C$ whenever the distortion is differentiable.
\begin{prop}[Pagès \cite{PAG1}]\label{prop:diff}
Under Assumption \ref{assum:compactdensity}, the distortion $C$ is continuously differentiable at every $w =(w_1, \ldots, w_\kappa) \in \Dcal_*^\kappa$.
Furthermore, for all $1 \leq \ell \leq \kappa$,
\[\nabla_\ell C(w) = \displaystyle \int_{W_\ell(w)}{ \left( w_\ell - \zbf \right) d\mu(\zbf) }.\]
\end{prop}

Some necessary conditions on the location of the minimizers of $C$ can be derived from its differentiability properties. Therefore, Proposition \ref{prop:strongMinima} below states that the minimizers of $C$ have parted components and that they are contained in the support of the density. Thus, the gradient is well defined and these minimizers are necessarily some zeroes of $\nabla C$. For the sequel it is convenient to let $\stackrel{\circ}{A}$ be the interior of any subset $A$ of $\rdkap$.\\
\begin{prop}[Pagès \cite{PAG1}]\label{prop:strongMinima}
Under Assumption \ref{assum:compactdensity}, we have
\begin{equation*}
\argmin_{w \in \left( \real ^d\right)^\kappa}C(w) \subset \argminloc_{w \in \Gcal^\kappa}C(w) \subset \; \stackrel{\circ}{\Gcal^\kappa} \cap \left\{\nabla C = 0\right\} \cap \Dcal_*^\kappa,
\end{equation*}
where $\argminloc_{w \in \Gcal^\kappa}C(w)$ stands for the set of local minimizers of $C$ over $\Gcal^\kappa$.
\end{prop}

For any $\zbf \in \real^d$ and $w \in \rdkap$, let us define the following vector of $\rdkap$
\begin{equation}\label{def:defH}
H(\zbf,w) \triangleq \left( \left(w_\ell - \zbf\right) \mathds{1}_{\left\{\zbf \in W_\ell(w)\right\}}\right)_{1 \leq \ell \leq \kappa}.
\end{equation}

On $\Dcal_*^\kappa$, the function $H$ may be interpreted as an observation of the gradient. With this notation, Proposition \ref{prop:diff} states that
\begin{equation}\label{assum:Hfunction}
\nabla C (w) = \displaystyle\int_\Gcal{ H(\zbf,w)d\mu(\zbf) }, \quad \mbox{$w \in \Dcal_*^\kappa$}.
\end{equation}

Let $\complement A$ stands for the complementary in $\rdkap$ of a subset $A \subset \rdkap$. Clearly, for all $w\in \complement \Dcal_*^\kappa$, the mapping $H(.,w)$ is integrable. Therefore, $\nabla C$ can be extended on $\rdkap$ \emph{via} the formula
\begin{equation}\label{eq:defh}
h(w) \triangleq \displaystyle\int_{\Gcal}{H(\zbf,w) d\mu(\zbf)}, \quad \mbox{$w \in \rdkap$}.
\end{equation}
Note however that the function $h$, which is sometimes called the average function of the algorithm, is not continuous.

\begin{rem}\label{rem:notLipchitz}
 Under Assumption \ref{assum:compactdensity}, a computation for all $w \in \Dcal_*^\kappa$ of the Hessian matrix $\nabla^2 C(w)$ can be deduced from Theorem 4 of Fort and Pagès \cite{FOR2}. In fact, the formula established in this theorem is valid for cost functions which are  more complex than $C$ (they are associated to Kohonen Self Organizing Maps, see Kohonen \cite{KOH1} for more details). In Theorem 4, letting $\sigma(k) = \mathds{1}_{\{k=0\}}$, provides the result for our distortion $C$. The resulting formula shows that $h$ is singular on $\complement\Dcal_*^\kappa$ and, consequently, that this function cannot be Lipschitz on $\Gcal^\kappa$.
\end{rem}

\subsection{Convergence of the CLVQ algorithm}\label{sec:VQ:Optimization}
The problem of finding a reliable clustering scheme for a dataset is equivalent to find optimal (or at least nearly optimal) minimizers for the mapping $C$. A minimization procedure by a usual gradient descent method cannot be implemented as long as $\nabla C$ is unknown. Thus, the gradient is approximated by a single example extracted from the data. This leads to the following stochastic gradient descent procedure
\begin{equation}\label{eq:gradientdescent}
w(t+1) = w(t) - \varepsilon_{t+1} H\left(\zbf_{t+1},w(t)\right), \quad \mbox{$t \geq 0$},
\end{equation}

where $w(0) \in \; \stackrel{\circ}{\Gcal^\kappa}\cap \; \Dcal_*^\kappa$ and $\zbf_1,\zbf_2 \ldots$ are independent observations distributed according to the probability measure $\mu$.\\

The algorithm defined by the iterations (\ref{eq:gradientdescent}) is known as the CLVQ algorithm in the data analysis community. It is also called the Kohonen Self Organizing Map algorithm with 0 neighbor (see for instance Kohonen \cite{KOH1}) or the on-line $k$-means procedure (see MacQueen \cite{MACQ} and Bottou \cite{BOT1}) in various fields related to statistics. As outlined by Pag\`es in \cite{PAG1}, this algorithm belongs to the class of stochastic gradient descent methods. However, the almost sure convergence of this type of algorithm cannot be obtained by general tools such as Robbins-Monro method (see Robbins and Monro \cite{ROB1}) or the Kushner-Clark's Theorem (see Kushner and Clark \cite{KUS1}). Indeed, the main difficulty essentially arises from the non convexity of the function $C$, its non coercivity and the singularity of $h$ at $\complement \Dcal_*^\kappa$ (we refer the reader to \cite[Section 6]{PAG1} for more details).\\

The following assumption set is standard in a gradient descent context. It basically upraises constraints on the decreasing speed of the sequence of steps $\left\{\varepsilon_t \right\}_{t=0}^{\infty}$.
\begin{assum}[Decreasing steps]\label{assum:deacreasingSteps}
 The $(0,1)$-valued sequence $\left\{\varepsilon_t \right\}_{t=0}^{\infty}$ satisfies the following two constraints:
\begin{enumerate}
\item $\sum_{t=0}^{\infty}{\varepsilon_t}  = \infty$.
\item $\sum_{t=0}^{\infty}{\varepsilon_t^2}  < \infty$.
\end{enumerate}
\end{assum}

\begin{figure}[!h]
\begin{center}
\includegraphics[scale = 0.50]{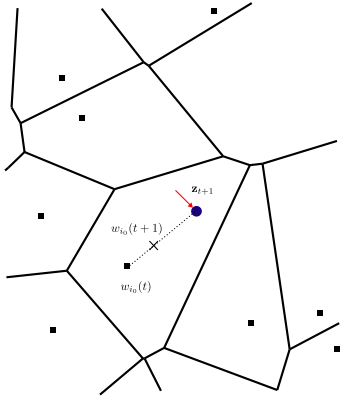}
\end{center}
\caption{Drawing of a portion of a 2-dimensional Voronoï tessellation. For $t\geq 0$, if $\zbf_{t+1} \in W_{\ell_0}\left(w(t)\right)$ then $w_\ell(t+1) = w_\ell(t)$ for all $\ell \neq \ell_0$ and $w_{\ell_0}(t+1)$ lies in the segment $\left[w_{\ell_0}(t), \zbf_{t+1}\right]$. The update of the vector $w_{\ell_0}(t)$ can also be viewed as a $\zbf_{t+1}$-centered homothety with ratio $1-\varepsilon_{t+1}$.}
\label{fig:homotethyCells}
\end{figure}

An examination of identities (\ref{eq:gradientdescent}) and (\ref{def:defH}) reveals that if $\zbf_{t+1} \in W_{\ell_0}\left(w(t)\right)$, where $\ell_0 \in \left\{ 1, \ldots, M\right\}$ then
\begin{equation}\label{eq:convexComb}
w_{\ell_0}(t+1) = \left(1 -\varepsilon_{t+1} \right)w_{\ell_0}(t) + \varepsilon_{t+1} \zbf_{t+1}.
\end{equation}
 The component $w_{\ell_0}(t+1)$ can be viewed as the image of $w_{\ell_0}(t)$ by a $\zbf_{t+1}$-centered homothety with ratio $1-\varepsilon_{t+1}$ (Figure \ref{fig:homotethyCells} provides an illustration of this fact). Thus, under Assumptions \ref{assum:compactdensity} and \ref{assum:deacreasingSteps}, the trajectories of $\left\{w(t)\right\}_{t=0}^{\infty}$ stay in $\stackrel{\circ}{\Gcal^\kappa}\cap \; \Dcal_*^\kappa$. More precisely,
  if \[w(0) \in \; \stackrel{\circ}{\Gcal^\kappa}\cap \; \Dcal_*^\kappa\]
  then
   \[w(t) \in \; \stackrel{\circ}{\Gcal^\kappa}\cap \; \Dcal_*^\kappa, \quad \mbox{$t \geq 0$, a.s.}\]

Although $\nabla C$ is not continuous some regularity can be obtained. To this end, we need to introduce the following materials. For any $\delta >0$ and any compact set $L \subset\real^d$, let the compact set $L^\kappa_\delta \subset \rdkap$ be defined as

\begin{equation}\label{eq:genCompact}
L^\kappa_\delta \triangleq \left\{w \in L^\kappa \; \vert \; \min_{k \neq \ell}\left\|w_\ell -w_k \right\| \geq \delta  \right\}.
\end{equation}
The next lemma that states on the regularity of $\nabla C$ will prove to be extremely useful in the proof of Theorem \ref{thm:pages} and throughout Section \ref{sect:DALVQ}.
\begin{lemma}[Pagès \cite{PAG1}]\label{lem:localLipschitz}
Assume that $\mu$ satisfies Assumption \ref{assum:compactdensity} and let $L$ be a compact set of $\real^d$. Then, there is some constant $P_\delta$ such that for all $w$ and $v$ in $L^\kappa_\delta$ with $[w, v] \subset \Dcal_*^\kappa$,
\[\left\| \nabla C(w) - \nabla C(v) \right\| \leq P_\delta \left\|w - v\right\|.\]
\end{lemma}

The following lemma, called G-Lemma in \cite{PAG1} is an easy-to-apply convergence results on stochastic algorithms. It is particularly adapted to the present situation of the CLVQ algorithm where the average function of the algorithm $h$ is singular.
\begin{theorem}[G-Lemma, Fort and Pagès \cite{FOR1}]\label{thm:GLemma}
Assume that the iterations (\ref{eq:gradientdescent}) of the CLVQ algorithm satisfy the following conditions:
\begin{enumerate}
\item $\sum_{t = 1}^{\infty}{\varepsilon_t} = \infty$ and $\varepsilon_t \xrightarrow[t \rightarrow \infty]{} 0$.
\item The sequences $\left\{w(t)\right\}_{t = 0}^{\infty}$ and $\left\{h\left(w(t)\right)\right\}_{t = 0}^{\infty}$ are bounded a.s.
\item The series $\sum_{t = 0}^{\infty}{\varepsilon_{t+1} \left(H(\zbf_{t+1},w(t)) - h(w(t)) \right)}$ converge a.s. in $\rdkap$.
    \item There exists a lower semi-continuous function $G:\rdkap \longrightarrow [0, \infty)$ such that
     \[\sum_{t =0}^{\infty}{\varepsilon_{t+1}G(w(t))} < \infty, \quad \mbox{a.s.}\]
\end{enumerate}
Then, there exists a random connected component $\Xi$ of $\left\{G=0 \right\}$ such that
\begin{equation*}
\dist\left(w(t),\Xi\right) \xrightarrow[t \rightarrow \infty]{} 0, \quad \mbox{a.s.,}
\end{equation*}
where the symbol $\dist$ denotes the usual distance function between a vector and a subset of $\rdkap$. Note also that if the connected components of $\left\{G=0 \right\}$ are singletons then there exists $\xi \in \left\{G=0 \right\}$ such that $w(t) \xrightarrow[t \rightarrow \infty]{} \xi$ a.s.
\end{theorem}
For a definition of the topological concept of connected component, we refer the reader to Choquet \cite{CHO2}. The interest of the G-Lemma depends upon the choice of $G$. In our context, a suitable lower semi-continuous function is $\widehat{G}$ defined by
\begin{equation}\label{eq:Ghat}
\widehat{G}(w) \triangleq \liminf_{v \in \Gcal^\kappa \cap \; \Dcal_*^\kappa, \; v \rightarrow w} \left\|\nabla C (v)\right\|^2, \qquad \mbox{$w \in \Gcal^\kappa$}.
\end{equation}

The next theorem is, as far as we know, the first almost sure convergence theorem for the stochastic algorithm CLVQ.
\begin{theorem}[Pag\`es \cite{PAG1}]\label{thm:pages}
Under Assumptions \ref{assum:compactdensity} and \ref{assum:deacreasingSteps}, conditioned on the event
\[\left\{\liminf_{t \rightarrow \infty}{\dist\left(w(t), \complement{\Dcal_*^\kappa }\right)} > 0 \right\}, \quad \mbox{one has}\]
\[\dist(w(t),\Xi_\infty) \xrightarrow[t \rightarrow \infty]{} 0, \quad \mbox{a.s.},\]
where $\Xi_\infty$ is some random connected component of $\left\{\nabla C = 0 \right\}$.
\end{theorem}
The proof is an application of the above G-Lemma with the mapping $\widehat{G}$ defined by equation (\ref{eq:Ghat}). Theorem \ref{thm:pages} states that the iterations of the CLVQ necessarily converge towards some critical points (zeroes of $\nabla C$). From Proposition \ref{prop:strongMinima} we deduce that the set of critical points necessarily contains optimal quantizers. Recall that without more assumption than $w(0) \in \; \stackrel{\circ}{\Gcal^\kappa}\cap \;\Dcal_*^\kappa$, we have already discussed the fact that the components of $w(t)$ are almost surely parted for all $t \geq 0$. Thus, it is easily seen that the two following events only differ on a set of zero probability
 \[\left\{\liminf_{t \rightarrow \infty}{\dist\left(w(t), \complement{\Dcal_*^\kappa }\right)} > 0 \right\}\]
 and
 \[\left\{\inf_{t \geq 0}{\dist\left(w(t), \complement{\Dcal_*^\kappa }\right)} > 0 \right\}.\]
Some results are provided by Pagès in \cite{PAG1} for asymptotically stuck components but, as pointed out by the author, they are less satisfactory.

\section{General distributed asynchronous algorithm}\label{sect:MDA}

\subsection{Model description}
Let $s(t)$ be any $\rdkap$-valued vector and consider the following iterations on a vector $w \in \rdkap$
\begin{equation}\label{eq:basiciteration}
w(t+1) = w(t) + s(t), \quad \mbox{$t\geq0$}.
\end{equation}

Here, the model of discrete time described by iterations (\ref{eq:basiciteration}) can only be performed by a single computing entity. Therefore, if the computations of the vectors $s(t)$ are relatively time consuming then not many iterations can be achieved for a given time span. Consequently, a parallelization of this computing scheme should be investigated. The aim of this section is to discuss a precise mathematical description of a distributed asynchronous model for the iterations (\ref{eq:basiciteration}). This model for distributed computing was originally proposed by Tsitsiklis et al. in \cite{TSI1} and was revisited in Bertsekas and Tsitsiklis \cite[Section 7.7]{BER1}.\\

Assume that we dispose of a distributed architecture with $M$ computing entities called processors (or agents, see for instance Bullo et al.~\cite{BUL1}). Each processor is labeled, for simplicity of notation, by a natural number $i \in \{1, \ldots, M\}$. Throughout the paper, we will add the superscript $i$ on the variables possessed by the processor $i$. In the model we have in mind, each processor has a buffer where its current version of the iterated vector is kept, i.e a local memory. Thus, for agent $i$ such iterations are represented by the $\rdkap$-valued sequence $\left\{w^i(t)\right\}_{t=0}^{\infty}$.\\

Let $t \geq 0$ denote the current time. For any pair of processors $(i,j) \in \left\{1, \ldots, M \right\}^2$, the value kept by agent $j$ and available for agent $i$ at time $t$ is not necessarily the most recent one, $w^j(t)$, but more probably and outdated one, $w^j(\tau^{i,j}(t))$, where the deterministic time instant $\tau^{i,j}(t)$ satisfy $0 \leq \tau^{i,j}(t) \leq t$. Thus, the difference $t - \tau^{i,j}(t)$ can be seen as a communication delay. This is a modeling of some aspects of the network: latency and bandwidth finiteness.\\
%

We insist on the fact that there is a distinction between ``global'' and ``local'' time. The time variable we refer above to as $t$ corresponds to a global clock. Such a global clock is needed only for analysis purposes. The processors work without knowledge of this global clock. They have access to a local clock or to no clock at all.\\

The algorithm is initialized at $t=0$, where each processor $i \in \left\{1, \ldots, M \right\}$ has an initial version $w^i(0) \in \rdkap$ in its buffer.
We define the general distributed asynchronous algorithm by the following iterations
\begin{equation}\label{eq:defwi}
w^i(t+1) =\sum_{j=1}^M a^{i,j}(t)w^j(\tau^{i,j}(t)) + s^i(t), \quad \mbox{$i \in \left\{1, \ldots, M \right\}$ and $t \geq 0$}.
\end{equation}

The model can be interpreted as follows: at time $t \geq 0$, processor $i$ receives messages from other processors containing $w^j(\tau^{i,j}(t))$. Processor $i$ incorporates these new vectors by forming a convex combination and incorporates the vector $s^i(t)$ resulting from its own ``local'' computations.
 The coefficients $a^{i,j}(t)$ are nonnegative numbers which satisfy the constraint
 \begin{equation}\label{eq:combiningcoeff}
 \sum_{j=1}^M a^{i,j}(t) =1, \quad \mbox{$i \in \left\{1, \ldots, M \right\}$ and $t \geq 0$}.
 \end{equation}

As the combining coefficients $a^{i,j}(t)$ depend on $t$, the network communication topology is sometimes referred to as time-varying. The sequences $\left\{\tau^{i,j}(t)\right\}_{t = 0}^\infty$ need not to be known in advance by any processor. In fact, their knowledge is not required to execute iterations defined by equation (\ref{eq:defwi}). Thus, we do not necessary dispose of a shared global clock or synchronized local clocks at the processors.\\

As for now the descent terms $\left\{s^i(t)\right\}_{t=0}^{\infty}$ will be arbitrary $\rdkap$-valued sequences. In Section \ref{sect:DALVQ}, when we define the Distributed Asynchronous Learning Vector Quantization (DALVQ), the definition of the descent terms will be made more explicit.\\

\begin{figure}[h]
\begin{center}
\includegraphics[scale = 0.35]{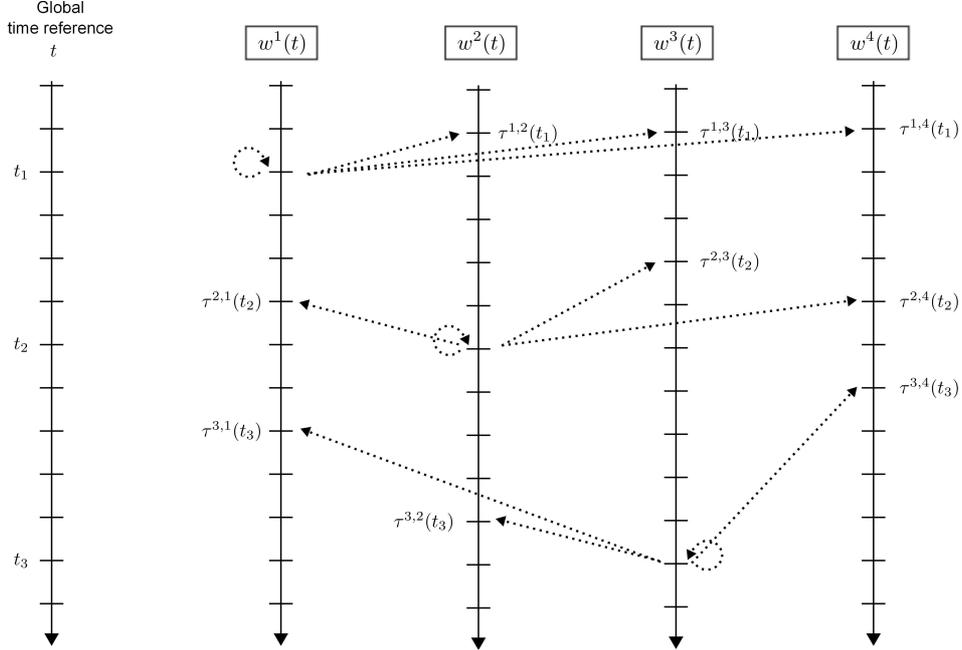}
\end{center}
\caption{Illustration of the time delays introduced in the general distributed asynchronous algorithm. Here, there are $M=4$ different processors with their own computations of the vectors $w^{(i)}$, $i \in \left\{1,2,3,4\right\}$. Three arbitrary values of the global time $t$ are represented ($t_1$, $t_2$ and $t_3$), with $\tau^{i,i}(t_k) = t_k$ for all $i \in \left\{1,2,3,4\right\}$ and $1 \leq k \leq 3$. The dashed arrows head towards the versions available at time $t_k$ for an agent $i \in \left\{1,2,3,4\right\}$ represented by the tail of the arrow.}
\label{fig:timescales}
\end{figure}

\subsection{The agreement algorithm}
This subsection is devoted to a short survey of the results, found by Blondel et al. in \cite{BLO1}, for a natural simplification of the general distributed asynchronous algorithm (\ref{eq:defwi}). This simplification is called agreement algorithm by Blondel et al. and is defined by
\begin{equation}\label{eq:agreementVector}
x^i(t+1) = \sum_{j=1}^M a^{i,j}(t)x^j(\tau^{i,j}(t)), \quad \mbox{$i \in \left\{1, \ldots, M \right\}$ and $t \geq 0$}.
\end{equation}
where $x^i(0) \in \rdkap$. An observation of these equations reveals that they are similar to iterations (\ref{eq:defwi}), the only difference being that all descent terms equal $0$.\\

In order to analyse the convergence of the agreement algorithm (\ref{eq:agreementVector}), Blondel et al in \cite{BLO1} define two sets of assumptions that enforce some weak properties on the communication delays and the network topology. As shown in \cite{BLO1}, if the assumptions contained in one of these two set hold, then the distributed versions of the agreement algorithm, namely the $x^i$, reach an asymptotical consensus. This latter statement means that there exists a vector $x^\star$ (independent of $i$) such that
\[ x^i(t) \xrightarrow[t \rightarrow \infty]{} x^\star, \quad \mbox{$i \in \left\{1, \ldots,M\right\}$.}\]

The agreement algorithm (\ref{eq:agreementVector}) is essentially driven by the communication times $\tau^{i,j}(t)$ assumed to be deterministic but do not need to be known \emph{a priori} by the processors. The following Assumption \ref{assum:boundedDelays} essentially ensures, in its third statement, that the communication delays $t -\tau^{i,j}(t)$ are bounded. This assumption prevents some processor from taking into account some arbitrarily old values computed by others processors. Assumption \ref{assum:boundedDelays} 1. is just a convention: when $a^{i,j}(t) = 0$ the value $\tau^{i,j}(t)$ has no effect on the update. Assumption \ref{assum:boundedDelays} 2. is rather natural because processors have access to their own most recent value.
\begin{assum}[Bounded communication delays]\label{assum:boundedDelays}
\begin{enumerate}
\item If $a^{i,j}(t) = 0$ then $\tau^{i,j}(t) =t$, \quad \mbox{$(i,j) \in \left\{1, \ldots, M \right\}^2$ and $t\geq 0$},
\item $\tau^{i,i}(t) = t, \quad i\in \left\{1, \ldots, M \right\} \text{ and }t\geq 0$.
\item There exists a positive integer $B_1$ such that
\[t-B_1 < \tau^{i,j}(t) \leq t, \quad \mbox{$(i,j) \in \left\{1, \ldots , M\right\}^2$ and $t \geq 0$.}\]
\end{enumerate}
\end{assum}

The next Assumption \ref{assum:stochMatrices} states that the value possessed by agent $i$ at time $t+1$, namely $x^i(t+1)$, is a weighted average of its own value and the values that it has just received from other agents.
\begin{assum}[Convex combination and threshold]\label{assum:stochMatrices}
There exists a positive constant $\alpha > 0$ such that the following three properties hold:
\begin{enumerate}
\item $a^{i,i}(t) \geq \alpha, \quad i \in \left\{1, \ldots, M\right\} \text { and }  t\geq 0$.
\item $a^{i,j}(t) \in \{0\} \cup [\alpha,1], \quad (i,j) \in \left\{1, \ldots, M \right\}^2 \text { and } t \geq 0$.
\item $\sum_{j=1}^M a^{i,j}(t) = 1, \quad i \in \left\{1, \ldots, M\right\} \text { and } t \geq 0$.
\end{enumerate}
\end{assum}

Let us mention one particular relevant case for the choice of the combining coefficients $a^{i,j}(t)$. Let $i \in \left\{1, \ldots, M\right\}$ and $t \geq 0$, the set \[N^i(t) \triangleq \left\{j \in \left\{1, \ldots,M \right\} \in \left\{1, \ldots, M\right\}  \; \vert \;  a^{i,j}(t) \neq 0 \right\}\]
 corresponds to the set of agents whose version is taken into account by processor $i$ at time $t$. For all $(i,j) \in \left\{1, \ldots, M \right\}^2$ and $t \geq 0$, the weights $a^{i,j}(t)$ are defined by
\[
a^{i,j}(t) =
\begin{cases}
1/\#N^i(t) & \text{ if } j \in N^i(t); \\
0 & \text{ otherwise;}
\end{cases}
\]
where $\#A$ denotes the cardinal of any finite set $A$. The above definition on the combining coefficients appears to be relevant for practical implementations of the model DALVQ introduced in Section \ref{sect:DALVQ}. For a discussion on others special interest cases regarding the choices of the coefficients $a^{i,j}(t)$ we refer the reader to \cite{BLO1}.\\

The communication patterns,  sometimes refereed to as the network communication topology, can be expressed in terms of directed graph. For a thorough introduction to graph theory, see Jungnickel \cite{JUN1}.
\begin{definition}[Communication graph]
Let us fix $t\geq 0$, the communication graph at time t, $\left(\mathcal{V}, E(t)\right)$, is defined by
\begin{itemize}
\item the set of vertices $\mathcal{V}$ is formed by the set of processors $\mathcal{V}= \left\{1, \ldots, M \right\}$,
\item the set of edges $E(t)$ is defined via the relationship
\[(j,i) \in E(t)\text{ if and only if }a^{i,j}(t) > 0.\]
\end{itemize}
\end{definition}

Assumption \ref{assum:GraphConnectivity} is a minimal condition required for a consensus among the processors. More precisely, it states that for any pair of agents $(i,j) \in \left\{1, \ldots,M \right\}^2$ there is a sequence of communications where the values computed by agent $i$ will influence (directly or indirectly) the future values kept by agent $j$.
\begin{assum}[Graph connectivity]\label{assum:GraphConnectivity}
The graph $\left(\mathcal{V}, \cup_{s \geq t}E(s)\right)$\\
is strongly connected for all $t\geq0$.
\end{assum}

Finally, we define two supplementary assumptions. The combination of one of the two following assumptions with the three previous ones will ensure the convergence of the agreement algorithm. As mentioned above, if Assumption \ref{assum:GraphConnectivity} holds then there is a communication path between any pair of agents. Assumption \ref{assum:BoundedInterCom} below expresses the fact that there is a finite upper bound for the length of such paths.
\begin{assum}[Bounded communication intervals]\label{assum:BoundedInterCom}
If $i$ communicates with $j$ an infinite number of times then there is a positive integer $B_2$ such that
\[(i,j) \in E(t) \cup E(t+1) \cup \ldots \cup E(t+ B_2 -1),\quad \mbox{$t \geq 0$.}\]
\end{assum}

Assumption \ref{assum:tminustau} is a symmetry condition: if agent $i \in \left\{1, \ldots, M\right\}$ communicates with agent $j \in \left\{1, \ldots, M\right\}$ then $j$ has communicated or will communicate with $i$ during the time interval $(t-B_3,t+B_3)$ where $B_3 > 0$.
\begin{assum}[Symmetry]\label{assum:tminustau}
There exists some $B_3 > 0$ such that whenever $(i,j) \in E(t)$, there exists some $\tau$ that satisfies $|t - \tau | < B_3$ and $(j,i) \in E(\tau)$.
\end{assum}

To shorten a little bit the notation, we set\\
\begin{tabular}{cc}
$
{ \bf \left(AsY\right)_1 } \equiv
\begin{cases}
\text{Assumption \ref{assum:boundedDelays};}\\
\text{Assumption \ref{assum:stochMatrices};}\\
\text{Assumption \ref{assum:GraphConnectivity};}\\
\text{Assumption \ref{assum:BoundedInterCom}.}\\
\end{cases}
$
&
$
{\bf \left(AsY\right)_2} \equiv
\begin{cases}
\text{Assumption \ref{assum:boundedDelays};}\\
\text{Assumption \ref{assum:stochMatrices};}\\
\text{Assumption \ref{assum:GraphConnectivity};}\\
\text{Assumption \ref{assum:tminustau};}\\
\end{cases}
$
\end{tabular}\\

We are now in a position to state the main result of this section. The Theorem \ref{thm:convAgreement} expresses the fact that, for the agreement algorithm, a consensus is asymptotically reached by the agents.
\begin{theorem}[Blondel et al.~\cite{BLO1}]\label{thm:convAgreement}
Under the set of Assumptions ${\bf \left(AsY\right)_1}$ or ${ \bf \left(AsY\right)_2}$,
there is a consensus vector $x^\star \in \rdkap$ (independent of $i$) such that
\[\lim_{t \rightarrow \infty}\left\|x^i(t) - x^\star \right\| = 0, \quad \mbox{$i \in \left\{1,\ldots,M\right\}$}.\]
Besides, there exist $\rho \in [0,1)$ and $L > 0$ such that
\[\left\|x^i(t) - x^i(\tau)\right\| \leq L\rho^{t - \tau}, \quad \mbox{$i \in \left\{1,\ldots,M\right\}$ and $t \geq \tau \geq 0$.} \]
\end{theorem}

\subsection{Asymptotic consensus}\label{sect:MDD}
This subsection is devoted to the analysis of the general distributed asynchronous algorithm (\ref{eq:defwi}). For this purpose, the study of the agreement algorithm defined by equations (\ref{eq:agreementVector}) will be extremely fruitful. The following lemma states that the version possessed by agent $i \in \left\{1, \ldots,M\right\}$ at time $t\geq 0$, namely $w^i(t)$, depends linearly on the others initialization vectors $w^j(0)$ and the descent subsequences $\left\{s^j(\tau)\right\}_{\tau=-1}^{t-1}$, where $j \in \left\{1, \ldots, M\right\}$.

\begin{lemma}[Tsitsiklis \cite{TSI2}]\label{lem:multDecomp}
For all  $(i,j) \in \left\{1, \ldots,M \right\}^2$ and $t \geq 0$, there exists a real-valued sequence $\left\{\phi^{i,j}\left(t,\tau\right)\right\}_{\tau = -1}^{t-1}$ such that
\[w^i(t) = \sum_{j=1}^M{\phi^{i,j}\left(t,-1\right)w^j(0)} + \sum_{\tau = 0}^{t-1}{\sum_{j=1}^M{\phi^{i,j}\left(t, \tau\right)s^j(\tau)}}.\]
\end{lemma}

For all $(i,j) \in \left\{1,\ldots,M\right\}^2$ and $t \geq 0$, the real-valued sequences $\left\{\phi^{i,j}\left(t,\tau\right)\right\}_{\tau = -1}^{t-1}$ do not depend on the value taken by the descent terms $s^i(t)$. The real numbers $\phi^{i,j}\left(t ,\tau \right)$ are determined by the sequences $\left\{\tau^{i,j}(\tau)\right\}_{\tau = 0}^t$ and $\left\{a^{i,j}(\tau)\right\}_{\tau=0}^t$ which do not depend on $w$. These last sequences are unknown in general, but some useful qualitative properties can be derived, as expressed in Lemma \ref{lem:decompProp} below.
\begin{lemma}[Tsitsiklis \cite{TSI2}]\label{lem:decompProp}
For all $(i,j) \in \left\{1, \ldots,M \right\}^2$, let $\left\{\phi^{i,j}\left(t,\tau\right)\right\}_{\tau = -1}^{t-1}$ be the sequences defined in Lemma \ref{lem:multDecomp}.
\begin{enumerate}
\item Under Assumption \ref{assum:stochMatrices},
\[0 \leq \phi^{i,j}\left(t, \tau\right) \leq 1, \quad \mbox{$(i,j) \in \left\{1, \ldots,M \right\}^2$ and $t > \tau \geq -1$}.\]
\item Under Assumptions ${ \bf \left(AsY\right)_1}$ or ${ \bf \left(AsY\right)_2}$, we have:
\begin{enumerate}
\item For all $(i,j) \in \left\{1, \ldots,M \right\}^2$ and $\tau \geq -1$, the limit of $\phi^{i,j}\left(t, \tau\right)$ as $t$ tends to infinity exists and is independent of $j$. It will be denoted $\phi^i(\tau)$.
\item There exists some $\eta > 0$ such that
\[\phi^i(\tau) > \eta, \quad \mbox{$i\in \left\{1, \ldots,M\right\}$ and $\tau \geq -1$}.\]
\item There exist a constant $A >0$ and $\rho \in (0,1)$ such that
\[\left|\phi^{i,j}\left(t, \tau\right) - \phi^i(\tau) \right| \leq A\rho^{t-\tau}, \quad \mbox{$(i,j)\in \left\{1, \ldots,M\right\}^2$ and $t > \tau \geq -1$}.\]
\end{enumerate}
\end{enumerate}
\end{lemma}

\begin{figure}[h]
\begin{center}
\includegraphics[scale=0.24]{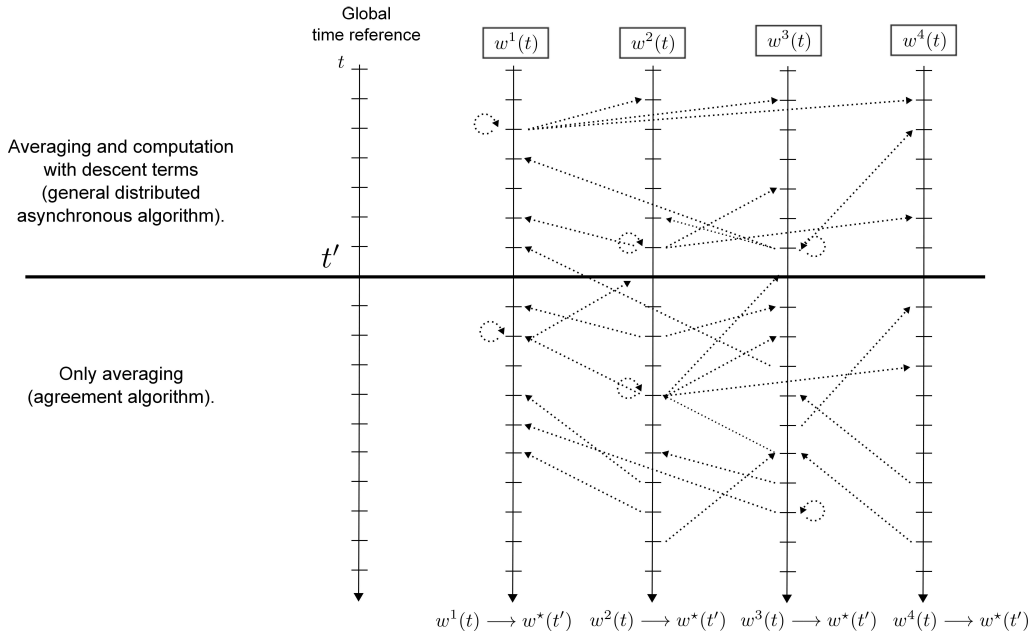}
\end{center}
\caption{The agreement vector at time $t'$, $w^\star(t')$ corresponds to the common value asymptotically achieved by all processors if computations integrating descent terms have stopped after $t'$, i.e, $s^j(t)=0$ for all $t \geq t'$.}
\label{fig:wstarAsConcensus}
\end{figure}

Take $t'\geq0$ and assume that the agents stop performing update after time $t'$, but keep communicating and merging the results. This means that $s^j(t)=0$ for all $t \geq t'$. Then, equations (\ref{eq:defwi}) write
\[w^i(t+1) = \sum_{j=1}^M{a^{i,j}(t)w^j\left(\tau^{i,j}(t)\right)}, \quad \mbox{$i \in \left\{1, \ldots, M \right\}$ and $t \geq t'$}.\]
If Assumptions ${\bf \left(AsY\right)_1}$ or ${\bf \left(AsY\right)_2}$ are satisfied then Theorem \ref{thm:convAgreement} shows that there is a consensus vector, depending on the time instant $t'$. This vector will be equal to $w^\star(t')$ defined below (see Figure \ref{fig:wstarAsConcensus}). Lemma \ref{lem:decompProp} provides a good way to define the sequence $\left\{w^\star(t)\right\}_{t=0}^\infty$ as shown in Definition \ref{eq:agreementVector}. Note that this definition does not involve any assumption on the descent terms.
\begin{definition}[Agreement vector]\label{def:agrementvector}
Assume that Assumptions ${\bf \left(AsY\right)_1}$ or ${ \bf \left(AsY\right)_2}$ are satisfied. The agreement vector sequence $\left\{w^\star(t)\right\}_{t=0}^{\infty}$ is defined by
\[w^\star(t) \triangleq \sum_{j=1}^M{\phi^j\left(-1\right)w^j(0)} + \sum_{\tau = 0}^{t-1}{\sum_{j=1}^M{\phi^j\left(\tau\right)s^j(\tau)}},\quad \mbox{$t \geq 0$}.\]
\end{definition}

It is noteworthy that the agreement vector sequence  $w^\star$ satisfies the following recursion formula
\begin{equation}\label{eq:agreemRecursion}
w^\star(t + 1) = w^\star(t) + \sum_{j=1}^M{\phi^j(t)s^j(t)}, \quad \mbox{$t \geq 0$}.
\end{equation}

\section{Distributed asynchronous learning vector quantization}\label{sect:DALVQ}
\subsection{Introduction, model presentation}
 From now on, and until the end of the paper, we assume that one of the two set of assumptions ${\bf \left(AsY\right)_1}$ or ${\bf \left(AsY\right)_2}$ holds, as well as the compact-supported density Assumption \ref{assum:compactdensity}. In addition, we will also assume that $0 \in \Gcal$. For the sake of clarity, all the proofs of the main theorems as well as the lemmas needed for these proofs have been postponed at the end of the paper, in Annex.\\

Tsitsiklis in \cite{TSI2}, Tsitsiklis et al in \cite{TSI1} and Bertsekas and Tsitsiklis in \cite{BER1} studied distributed asynchronous stochastic gradient optimization algorithms. In this series of publications, for the distributed minimization of a cost function $F : \rdkap \longrightarrow \real$, the authors considered the general distributed asynchronous algorithm defined by equation (\ref{eq:defwi}) with specific choices for stochastic descent terms $s^i$. Using the notation of Section \ref{sect:MDA}, the algorithm writes
\[w^i(t+1) = \sum_{j=1}^M{a^{i,j}(t)w^j(\tau^{i,j}(t))} +s^i(t), \quad \mbox{$i \in \left\{1, \ldots, M \right\}$ and $t \geq 0$,}\]
with stochastic descent terms $s^i(t)$ satisfying

\begin{align}\label{eq:tsitsiklischoice}
&\esp\left\{ s^i(t) \; {\big \vert} \; s^j(\tau), \; j \in \left\{1, \ldots,M \right\} \text{ and } t > \tau \geq 0\right\} = -\varepsilon_{t+1}^i \nabla F\left(w^i(t)\right), \notag\\
& \qquad \qquad \qquad \qquad \qquad \qquad \qquad \qquad \mbox{$i \in \left\{1, \ldots, M \right\}$ and $t \geq 0$}.
\end{align}
where $\left\{\varepsilon_{t}^i\right\}_{t=0}^{\infty}$ are decreasing steps sequences. The definition of the descent terms in \cite{BER1,TSI1} is more general than the one appearing in equation (\ref{eq:tsitsiklischoice}). We refer the reader to Assumption 3.2 and 3.3 in \cite{TSI1} and Assumption 8.2 in \cite{BER1} for the precise definition of the descent terms in \cite{BER1,TSI1}. As discussed in Section \ref{sect:CLVQ}, the CLVQ algorithm is also a stochastic gradient descent procedure. Unfortunately, the results from Tsitisklis et al. in \cite{TSI1,TSI2,BER1} do not apply with our distortion function, $C$, since the authors assume that $F$ is continuously differentiable and $\nabla F$ is Lipschitz. Therefore, the aim of this section is to extend the results of Tsitsiklis et al. to the context of vector quantization and on-line clustering. \\

We first introduce the Distributed Asynchronous Learning Vector Quantization (DALVQ) algorithm. To prove its almost sure consistency, we will need an Asynchronous G-Lemma, which is inspired from the G-Lemma, Theorem \ref{thm:GLemma}, presented in Section \ref{sect:CLVQ}. This theorem may be seen as an easy-to-apply tool for the almost sure consistency of a distributed asynchronous system where the average function is not necessary regular. Our approach sheds also some new light on the convergence of distributed asynchronous stochastic gradient descent algorithms. Precisely, Proposition 8.1 in \cite{TSI1} claims that $\liminf_{t \rightarrow \infty} \left\| \nabla F(w^i(t)) \right\| = 0$ while our main Theorem \ref{thm:AsyncTheorem} below states that $\lim_{t \rightarrow \infty} \left\| \nabla C(w^i(t)) \right\| = 0$. However, there is a price to pay for this more precise result with the non Lipschitz gradient $\nabla C$. Similarly to Pagès \cite{PAG1}, who assumes that the trajectory of the CLVQ algorithm has almost surely asymptotically parted components (see Theorem \ref{thm:pages} in Section \ref{sect:CLVQ}), we will suppose that the agreement vector sequence has, almost surely, asymptotically parted component trajectories.\\

Recall that the goal of the DALVQ is to provide a well designed distributed algorithm that processes quickly (in term of wall clock time) very large data sets to produce accurate quantization. The data sets (or streams of data) are distributed among several queues sending data to the different processors of our distributed framework. Thus, in this context the sequence $\zbf_{1}^i, \zbf_{2}^i, \ldots$ stands for the data available for processor, where $i \in \left\{1, \ldots, M \right\}$. The random variables
\[\zbf_{1}^1,\zbf_{2}^1, \ldots ,\zbf_1^2, \zbf_2^2, \ldots \]
are assumed to be independent and identically distributed according to $\mu$.\\

In the definition of the CLVQ procedure (\ref{eq:gradientdescent}), the term $H\left(\zbf_{t+1},w(t)\right)$ can be seen as an observation of the gradient $\nabla C\left(w(t)\right)$. Therefore, in our DALVQ algorithm, each processor $i \in \left\{1, \ldots, M \right\}$ is able to compute such observations using its own data $\zbf_1^i, \zbf_2^i, \ldots $. Thus, the DALVQ procedure is defined by equation (\ref{eq:defwi}) with the following choice for the descent term $s^i$:
\begin{equation}\label{eq:defSiDescent}
s^i(t) =
\begin{cases}
-\varepsilon^i_{t +1} H\left(\zbf_{t+1}^i, w^i(t) \right) & \text{ if } t \in T^i; \\
0 & \text{ otherwise;}
\end{cases}
\end{equation}
where $\left\{\varepsilon_{t}^i\right\}_{t=0}^{\infty}$ are $(0,1)$-valued sequences. The sets $T^i$ contain the time instants where the version $w^i$, kept by processor $i$, is updated with the descent terms. This fine grain description of the algorithm allows some processors to be idle for computing descent terms (when $t \notin T^i$). This reflects the fact that the computing operations might not take the same time for all processors, which is precisely the core of asynchronous algorithms analysis. Similarly to time delays and combining coefficients, the sets $T^i$ are supposed to be deterministic but do not need to be known \emph{a priori} for the execution of the algorithm.\\

In the DALVQ model, randomness arises from the data $\zbf$. Therefore, it is natural to let $\left\{\mathcal{F}_t\right\}_{t=0}^{\infty}$ be the filtration built on the $\sigma$-algebras
\[\mathcal{F}_t \triangleq \sigma\left(\zbf^i_s,\; i \in \left\{1, \ldots, M \right\} \text{ and } t \geq s \geq 0 \right), \quad \mbox{$t\geq 0$.}\]

An easy verification shows that, for all $j \in \left\{1, \ldots, M\right\}$ and $t\geq0$, $w^\star(t)$ and $w^j(t)$ are $\mathcal{F}_t$-measurable random variables.\\

For simplicity, the assumption on the decreasing speed of the sequences $\left\{\varepsilon^i_t\right\}_{t=0}^{\infty}$ is strengthened as follows. The notation $a \vee b$ stands for the maximum of two reals $a$ and $b$.
\begin{assum}\label{assum:strongdeacreasingSteps}
There exist two real numbers $K_1 > 0$ and $K_2 \geq 1$ such that
\[\frac{K_1}{t \vee 1} \leq \varepsilon^i_{t+1} \leq \frac{K_2}{t \vee 1}, \quad \mbox{$i \in \left\{1, \ldots, M \right\}$ and $t \geq 0$}.\]
\end{assum}
If Assumption \ref{assum:strongdeacreasingSteps} holds then the sequences $\left\{\varepsilon_t^i\right\}_{t=0}^{\infty}$ satisfy the standard Assumption \ref{assum:deacreasingSteps} for stochastic optimization algorithms. Note that the choice of steps proportional to $1/t$ has been proved to be a satisfactory learning rate, theoretically speaking and also for practical implementations (see for instance Murata \cite{MUR1} and Bottou and LeCun \cite{BOT5}).\\

For practical implementation, the sequences $\left\{\varepsilon_{t+1}^i\right\}_{t=0}^{\infty}$ satisfying Assumption \ref{assum:strongdeacreasingSteps} can be implemented without a global clock, that is, without assuming that the current value of $t$ is known by the agents. This assumption is satisfied, for example, by taking the current value of $\varepsilon_t^i$ proportional to $1/n^i_t$, where $n^i_t$ is the number of times that processor $i$ as performed an update, i.e., the cardinal of the set $T^i \cap \left\{ 0, \ldots, t\right\}$. For a given processor, if the time span between consecutive updates is bounded from above and from below, a straightforward examination shows that the sequence of steps satisfy Assumption \ref{assum:strongdeacreasingSteps}.\\

Finally, the next assumption is essentially technical in nature. It enables to avoid time instants where all processors are idle. It basically requires that, at any time $t \geq0$, there is at least one processor $i \in \left\{1, \ldots, M \right\}$ satisfying $s^i(t) \neq 0$.
\begin{assum}\label{assum:titimes}
One has $\sum_{j =1}^M{\Ind_{\left\{t \in T^j\right\}}} \geq1$ for all $t\geq0$.
\end{assum}

\subsection{The asynchronous G-Lemma}
The aim of this subsection is to state a useful theorem similar to Theorem \ref{thm:GLemma}, but adapted to our asynchronous distributed context. The precise Definition \ref{def:agrementvector} of the agreement vector sequence should not cast aside the intuitive definition. The reader should keep in mind that the vector $w^\star(t)$ is also the asymptotical consensus if descent terms are zero after time $t$. Consequently, even if the agreement vector $\left\{w^\star(t)\right\}_{t=0}^{\infty}$ is adapted to the filtration $\left\{\mathcal{F}_t\right\}_{t = 0}^{\infty}$, the vector $w^\star(t)$ cannot be accessible for a user at time $t$. Nevertheless, the agreement vector $w^\star(t)$ can be interpreted as a ``probabilistic state'' of the whole distributed quantization scheme at time $t$. This explains why the agreement vector is a such convenient tool for the analysis of the DALVQ convergence and will be central in our adaptation of G-Lemma, Theorem \ref{thm:AGLemma}.\\

Let us remark that equation (\ref{eq:agreemRecursion}), writes for all $t\geq0$,
\begin{align*}
w^\star(t+1) &= w^\star(t) + \sum_{j=1}^M{\phi^j(t)s^j(t)} \\
&=w^\star(t) - \sum_{j=1}^M\mathds{1}_{\left\{t \in T^j\right\}}\phi^j(t)\varepsilon_{t+1}^jH\left(\zbf_{t+1}^j,w^j(t)\right).
\end{align*}
We recall the reader that the $[0,1]$-valued functions $\phi^j$'s are defined in Lemma \ref{lem:multDecomp}.\\

Using the function $h$ defined by identity (\ref{eq:defh}) and the fact that the random variables $w^\star(t)$ and $w^j(t)$ are $\mathcal{F}_t$-measurable then it holds
\[h(w^\star(t)) = \esp\left\{H\left(\zbf, w^\star(t) \right)  \; \vert \; \mathcal{F}_t \right\}, \quad \mbox{$t \geq0$.}\]
and
\[h(w^j(t)) = \esp\left\{H\left(\zbf, w^j(t) \right) \; \vert \; \mathcal{F}_t \right\}, \quad \mbox{$j \in\left\{1, \ldots,M\right\}$ and $t\geq0$}.\]
where $\zbf$ is a random variable of law $\mu$ independent of $\mathcal{F}_t$.

For all $t \geq 0$, set
\begin{equation}\label{eq:defepsilonstar}
\varepsilon^\star_{t+1} \triangleq \sum_{j=1}^M{\Ind_{\left\{t \in T^j\right\}} \phi^j(t)\varepsilon_{t+1}^j}.
\end{equation}

Clearly, the real numbers $\varepsilon^\star_{t}$ are nonnegative. Their strictly positiveness will be discussed in Proposition \ref{prop:epsilonstarspeed}.

Set
\begin{equation}\label{eq:seriesM1}
\Delta M^{(1)}_t \triangleq \sum_{j=1}^M{\Ind_{\left\{t \in T^j\right\}}\phi^j(t) \varepsilon_{t+1}^j \left(h(w^\star(t)) - h(w^j(t))\right)}, \quad \mbox{$t \geq0$,}
\end{equation}
and
\begin{equation}\label{eq:seriesM2}
\Delta M^{(2)}_t \triangleq \sum_{j=1}^M{\Ind_{\left\{t \in T^j\right\}}\phi^j(t) \varepsilon_{t+1}^j\left(h(w^j(t)) -  H\left(\zbf_{t+1}^j,w^j(t) \right)\right)}, \quad \mbox{$t \geq0$.}
\end{equation}

Note that $\esp\left\{\Delta M^{(2)}_t \right\}=0$ and, consequently, that the random variables $\Delta M^{(2)}_t$ can be seen as the increments of a martingale with respect to the filtration $\left\{\mathcal{F}_t\right\}_{t=0}^{\infty}$.\\

Finally, with this notation, equation (\ref{eq:agreemRecursion}) takes the form
\begin{equation}\label{eq:agreemRecursionPages}
w^\star(t+1) = w^\star(t) - \varepsilon_{t+1}^\star h(w^\star(t)) + \Delta M^{(1)}_t + \Delta M^{(2)}_t, \quad \mbox{$t \geq 0$.}
\end{equation}

We are now in a position to state our most useful tool, which is similar in spirit to the G-Lemma, but adapted to the context of distributed asynchronous stochastic gradient descent algorithm.
\begin{theorem}[Asynchronous G-Lemma]\label{thm:AGLemma}
Assume that ${\bf \left(AsY\right)_1}$ or ${\bf \left(AsY\right)_2}$ and Assumption \ref{assum:compactdensity} hold and that the following conditions are satisfied:
\begin{enumerate}
\item$\sum_{t = 0}^{\infty}{\varepsilon^\star_t } = \infty$ and  $\varepsilon^\star_t \xrightarrow[t \rightarrow \infty]{} 0$.
\item The sequences $\left\{w^\star(t)\right\}_{t=0}^{\infty}$ and $\left\{h(w^\star(t))\right\}_{t=0}^{\infty}$ are bounded a.s.
\item The series $\sum_{t=0}^{\infty}{\Delta M_t^{(1)}}$ and  $\sum_{t=0}^{\infty}{\Delta M_t^{(2)}}$ converge a.s. in $\rdkap$.
\item There exists a lower semi-continuous function $G:\rdkap \longrightarrow [0, \infty)$ such that
\[\sum_{t = 0}^{\infty}{\varepsilon_{t+1}^\star}G\left(w^\star(t)\right) < \infty, \quad \mbox{a.s.}\]
\end{enumerate}
Then, there exists a random connected component $\Xi$  of $\left\{G=0 \right\}$ such that
\[\dist\left(w^\star(t),\Xi\right) \xrightarrow[t \rightarrow \infty]{} 0, \quad \mbox{a.s.}\]
\end{theorem}

\subsection{Trajectory analysis}
The Pag\`es's proof in \cite{PAG1} on the almost sure convergence of the CLVQ procedure required a careful examination of the trajectories of the process $\left\{w(t)\right\}_{t=0}^\infty$. Thus, in this subsection we investigate similar properties and introduce the assumptions that will be needed to prove our main convergence result, Theorem \ref{thm:AsyncTheorem}.

The next Assumption \ref{assum:witrajectories} ensures that, for each processor, the quantizers stay in the support of the density.
\begin{assum}\label{assum:witrajectories}
One has
\[\mathbb{P}\left\{w^j(t) \in \Gcal^\kappa \right\}=1, \quad \mbox{ $j \in \left\{1, \ldots, M \right\}$ and $t\geq 0$.}\]
\end{assum}

Firstly, let us mention that since the set $\Gcal^\kappa$ is convex, if Assumption \ref{assum:witrajectories} holds then
\[\mathbb{P}\left\{w^\star(t) \in \Gcal^\kappa \right\}=1, \quad \mbox{$t\geq0$}.\]

Secondly, note that the Assumption \ref{assum:witrajectories} is not particularly restrictive. This assumption is satisfied under the condition: for each processor, no descent term is added while a combining computation is performed. This writes
\[a_{i,j}(t) = \delta_{i,j} \text{ and }\tau^{i,i}(t)= t, \quad \mbox{ $(i,j) \in \left\{ 1, \ldots, M \right\}^2$ and $t \in T^i$}. \]
This requirement makes sense for practical implementations.

Recall that if $t \notin T^i$, then $s^i(t) =0$. Thus, equation (\ref{eq:defwi}) takes the form

\begin{equation}
w^i(t+1) =
\begin{cases}
\begin{array}{cc}
w^i(t+1) &= w^i(t) - \varepsilon_{t+1}^i\left( w^i(t) - \zbf^i_{t+1} \right) \\
&= \left(1 - \varepsilon_{t+1}^i \right)w^i(t) + \varepsilon_{t+1}^i \zbf^i_{t+1}
\end{array} & \text{ if } t \in T^i; \\
\begin{array}{cc}
w^i(t+1) &= \sum_{j=1}^M a^{i,j}(t)w^j(\tau^{i,j}(t))
\end{array}& \text{ otherwise.}
\end{cases}
\end{equation}

Since $\Gcal^\kappa$ is a convex set, it follows easily that if $w^j(0) \in \Gcal^\kappa$, then  $w^j(t) \in \Gcal^\kappa$ for all $j \in \left\{1, \ldots, M \right\}$ and $t\geq 0$ and, consequently, that Assumption \ref{assum:witrajectories} holds.\\

The next Lemma \ref{lem:normMaj} provides a deterministic upper bound on the differences between the distributed versions $w^i$ and the agreement vector. For any subset $A$ of $\rdkap$, the notation $\diam(A)$ stands for the usual diameter defined by
\[\diam(A) = \sup_{x,y \in A} \left\{ \left\|x - y\right\| \right\}.\]

\begin{lemma}\label{lem:normMaj}
Assume ${\bf \left(AsY\right)_1}$ or ${\bf \left(AsY\right)_2}$ holds and that Assumptions \ref{assum:compactdensity}, \ref{assum:strongdeacreasingSteps} and \ref{assum:witrajectories} are satisfied then
\[\|w^\star(t) - w^i(t)\| \leq \sqrt{\kappa}M\diam(\Gcal)A K_2\theta_t, \quad \mbox{$i \in \left\{1, \ldots, M \right\}$ and $t\geq 0$, a.s.,}\]
where $\theta_t \triangleq  \sum_{\tau=-1}^{t-1}\frac{1}{\tau \vee 1} \rho^{t- \tau}$, $A$ and $\rho$ are the constants introduced in Lemma \ref{lem:decompProp}, $K_2$ is defined in Assumption \ref{assum:strongdeacreasingSteps}.
\end{lemma}

The sequence $\left\{\theta_t\right\}_{t=0}^{\infty}$ defined in Lemma \ref{lem:normMaj} satisfies

\begin{equation}\label{eq:thetaprop}
\theta_t \xrightarrow[t\rightarrow \infty]{} 0  \text{ and } \sum_{t=0}^{\infty}\frac{\theta_t}{t} < \infty.
\end{equation}
We give some calculations justifying the statements at the end of the Annex.

Thus, under Assumptions \ref{assum:strongdeacreasingSteps} and \ref{assum:witrajectories}, it follows easily that
\[w^\star(t) -w^i(t) \xrightarrow[t \rightarrow \infty]{} 0, \quad \mbox{$i \in \left\{1, \ldots,M\right\}$, a.s.,}\]
and
\begin{equation}\label{eq:convwjwi}
w^i(t) -w^j(t) \xrightarrow[t \rightarrow \infty]{} 0, \quad \mbox{$(i,j) \in \left\{ 1, \ldots, M\right\}^2$, a.s.}
\end{equation}

This shows that the trajectories of the distributed versions of the quantizers reach asymptotically a consensus with probability 1. In other words, if one of the sequences $\left\{w^i(t)\right\}_{t=0}^{\infty}$ converges then they all converge towards the same value.
The rest of the paper is devoted to prove that this common value is in fact a zero of $\nabla C$, i.e. a critical point.\\

 To prove the result mentioned above, we will need the following assumption, which basically states that the components of $w^\star$ are parted, for every time $t$ but also asymptotically. This assumption is similar in spirit to the main requirement of Theorem \ref{thm:pages}.

\begin{assum}\label{assum:wstartrajectories}
One has
\begin{enumerate}
\item $\mathbb{P}\left\{w^\star(t) \in \Dcal_*^\kappa \right\}=1, \quad t\geq 0$.
\item $\mathbb{P}\left\{\liminf_{t \rightarrow \infty} \dist\left( w^\star(t), \complement \Dcal_*^\kappa \right) > 0 \right\}=1, \quad t\geq 0$.
\end{enumerate}
\end{assum}

\subsection{Consistency of the DALVQ}
In this subsection we state our main theorem on the consistency of the DALVQ. Its proof is based on the Asynchronous G-Lemma, Theorem \ref{thm:AGLemma}. The goal of the next proposition is to ensure that the first assumption of Theorem \ref{thm:AGLemma} holds.

\begin{prop}\label{prop:epsilonstarspeed}
Assume ${\bf \left(AsY\right)_1}$ or ${\bf \left(AsY\right)_2}$ holds and that Assumptions \ref{assum:compactdensity}, \ref{assum:strongdeacreasingSteps} and \ref{assum:titimes} are satisfied then $\varepsilon^\star_t > 0$, $t \geq 0$, $\varepsilon^\star_t \xrightarrow[t \rightarrow \infty]{} 0$ and $\sum_{t = 0}^{\infty}{\varepsilon^\star_t } = \infty$.
\end{prop}

The second condition required in Theorem \ref{thm:AGLemma} deals with the convergence of the two series defined by equations (\ref{eq:seriesM1}) and (\ref{eq:seriesM2}). The next Proposition \ref{prop:assumptionImplcations1} provides sufficient condition for the almost sure convergence of these series.
\begin{prop}\label{prop:assumptionImplcations1}
Assume ${\bf \left(AsY\right)_1}$ or ${\bf \left(AsY\right)_2}$ holds and that Assumptions \ref{assum:compactdensity}, \ref{assum:strongdeacreasingSteps}, \ref{assum:witrajectories} and \ref{assum:wstartrajectories} are satisfied then the series\\
$\sum_{t=0}^{\infty}{\Delta M_t^{(1)}}$ and  $\sum_{t=0}^{\infty}{\Delta M_t^{(2)}}$ converge almost surely in $\rdkap$.
\end{prop}

This next proposition may be considered has the most important step in the proof of the convergence of the DALVQ. It establishes the convergence of a series of the form $\sum_{t = 0}^{\infty}{\varepsilon_{t+1}\left\|\nabla C \left(w(t)\right)\right\|^2}$. The analysis of the convergence of this type of series is standard for the analysis of stochastic gradient method (see for instance Benveniste et al.~\cite{BEN1} and Bottou \cite{BOT2}). In our context, we pursue the fruitful use of the agreement vector sequence, $\left\{w^\star(t) \right\}_{t=0}^\infty$, and its related ``steps'', $\left\{\varepsilon_t^\star\right\}_{t=0}^{\infty}$.\\

Note that under Assumption \ref{assum:wstartrajectories}, we have $h\left(w^\star(t)\right) = \nabla C \left(w^\star(t)\right)$ for all $t \geq 0$, almost surely, therefore the sequence $\left\{\nabla C \left(w^\star(t)\right)\right\}_{t=0}^{\infty}$ below is well defined.
\begin{prop}\label{prop:convSeries}
Assume ${\bf \left(AsY\right)_1}$ or ${\bf \left(AsY\right)_2}$ holds and that Assumptions \ref{assum:compactdensity}, \ref{assum:strongdeacreasingSteps}, \ref{assum:witrajectories} and \ref{assum:wstartrajectories} are satisfied then
\begin{enumerate}
\item $C\left(w^\star(t)\right) \xrightarrow[ t \rightarrow \infty]{} C_\infty,\quad \mbox{a.s.,}$\\ where $C_\infty$ is a $[0, \infty)$-valued random variable,
\item
\begin{equation}\label{eq:statement}
\sum_{t=0}^{\infty}{\varepsilon_{t+1}^\star\left\|\nabla C \left(w^\star(t)\right)\right\|^2} < \infty, \quad \mbox{a.s.}
\end{equation}
\end{enumerate}
\end{prop}
 Remark that from the convergence of the series given by equation (\ref{eq:statement}) one can only deduce that $\liminf_{t \rightarrow \infty} \left\|\nabla C \left(w^\star(t) \right)\right\| = 0$.\\

We are now in a position to state the main theorem of this paper, which expresses the convergence of the distributed version towards some zero of the gradient of the distortion. In addition, the convergence results (\ref{eq:convwjwi}) imply that if a version converges then all the versions converge towards this value.
\begin{theorem}[Asynchronous Theorem]\label{thm:AsyncTheorem}
Assume ${\bf \left(AsY\right)_1}$ or ${\bf \left(AsY\right)_2}$ holds and that Assumptions \ref{assum:compactdensity}, \ref{assum:strongdeacreasingSteps}, \ref{assum:titimes}, \ref{assum:witrajectories} and \ref{assum:wstartrajectories} are satisfied
then
\begin{enumerate}
\item $w^*(t) -w^i(t) \xrightarrow[t \rightarrow \infty]{} 0, \quad \mbox{$i \in \left\{ 1, \ldots, M\right\}$, a.s.,}$
\item $w^i(t) -w^j(t) \xrightarrow[t \rightarrow \infty]{} 0, \quad \mbox{$(i,j) \in \left\{ 1, \ldots, M\right\}^2$, a.s.,}$
\item $\dist\left(w^\star(t),\Xi_\infty\right) \xrightarrow[t \rightarrow \infty]{} 0, \quad \mbox{a.s.,}$
\item $\dist\left(w^i,\Xi_\infty\right) \xrightarrow[t \rightarrow \infty]{} 0, \quad \mbox{$i \in \left\{1, \ldots, M\right\}$, a.s.,}$
\end{enumerate}
where $\Xi_\infty$ is some random connected component of the set $\left\{\nabla C = 0\right\} \cap \Gcal^\kappa$.
\end{theorem}

\subsection{Annex}
{\bf Sketch of the proof of Asynchronous G-Lemma \ref{thm:AGLemma}.} The proof is an adaptation of the one found by Fort and Pagès, Theorem 4 in \cite{FOR1}. The recursive equation (\ref{eq:agreemRecursionPages}) satisfied by the sequence $\left\{w^\star(t) \right\}_{t=0}^\infty$ is similar to the iterations (2) in \cite{FOR1}
(with the notation of this paper):
\[X^{t+1} = X^t - \varepsilon_{t+1}h\left(X^t\right) + \varepsilon_{t+1}\left(\Delta M^{t+1} + \eta^{t+1}\right), \quad \mbox{$t\geq0$}.\]

Thus, similarly, we define a family of continuous time stepwise function $\left\{u \mapsto \wch\left(t,u\right) \right\}_{t = 1}^\infty$.
\[\wch^\star\left(0,u\right) \triangleq w^\star(s)\text{, if } u \in [\varepsilon^\star_{1} + \ldots + \varepsilon^\star_s, \varepsilon^\star_{1} + \ldots + \varepsilon^\star_{s+1}), \quad \mbox{$u \in [0,\infty)$}.\]
and if $u < \varepsilon^\star_{1}$, $\wch^\star\left(0,u\right) = w^\star(0)$.
\[\wch^\star\left(t,u\right) \triangleq \wch^\star\left(0,\varepsilon^\star_{1} + \ldots + \varepsilon^\star_t + u\right), \quad \mbox{$t \geq 1$ and $u \in [0,\infty)$}.\]

Hence, for every $t \in \mathbb{N}$,
\[\wch^\star(t,u) = \wch^\star(0,t) - \displaystyle\int_{0}^{u}{h\left(\wch^\star(t,v)\right)dv} + R_u(t), \quad \mbox{$u \in [0,\infty)$},\]
where, for every $t \geq 1$ and $u \in [\varepsilon^\star_{1} + \ldots + \varepsilon^\star_{t + t'}, \varepsilon^\star_{1} + \ldots + \varepsilon^\star_{t + t' + 1} )$,
\[R_u(t) \triangleq \displaystyle\int_{\varepsilon^\star_{t} + \ldots + \varepsilon^\star_{t + t'}}^{\varepsilon^\star_{1} + \ldots + \varepsilon^\star_{t} + u}{\wch^\star(0,v)dv} + \sum_{s=t+1}^{t+t'}{\left(\Delta M^{(1)}_s + \Delta M^{(2)}_s\right)}.\]

The only difference between the families of continuous time functions $\left\{\wch\left(t,u\right) \right\}_{t = 1}^\infty$ and $\left\{X^{(t)}\right\}_{t=1}^\infty$ defined in \cite{FOR1} is the remainder term $R_u(t)$.
The convergence
\[\sup_{u \in [0,T]}\left\|R_u(t)\right\| \xrightarrow[t \rightarrow \infty]{} 0, \quad \mbox{$T >0$}.\] follows easily from the third assumption of Theorem \ref{thm:AGLemma}. The rest of the proof follows similarly as in Theorem 4 \cite{FOR1}.
\proofend

{\bf Proof of Lemma \ref{lem:normMaj}.} For all $i \in \left\{1, \ldots, M \right\}$, and all $t \geq 0$, and all $1 \leq \ell \leq \kappa$, we may write
\begin{align*}
&\left\|w_\ell^i(t) - w_\ell^\star(t)\right\| \\
& \quad =\left\| \sum_{j=1}^M \left(\left(\phi^{i,j}(t,-1) - \phi^j(-1) \right) w_\ell^j(0) + \sum_{\tau = 0}^{t-1} \left(\phi^{i,j}(t,\tau) - \phi^j(t) \right) s_\ell^j(\tau) \right)\right\|\\
& \quad \qquad (\mbox{by Definition \ref{def:agrementvector} and Lemma \ref{lem:multDecomp}})\\
& \quad \leq \sum_{j=1}^M \left|\phi^{i,j}(t,-1) - \phi^j(-1) \right|  \left\| w_\ell^j(0) \right\| +  \sum_{\tau = 0}^{t-1}  \sum_{j=1}^M \left|\phi^{i,j}(t,\tau) - \phi^j(t) \right| \left\|s_\ell^j(\tau) \right\| \\
& \quad \leq A \rho^{t+1}\sum_{j=1}^M  \left\| w_\ell^j(0) \right\| +   A \sum_{\tau = 0}^{t-1}  \sum_{j=1}^M \rho^{t-\tau} \left\|s_\ell^j(\tau) \right\|\\
& \quad \qquad (\mbox{by Lemma \ref{lem:decompProp}}).
\end{align*}
Thus,
\begin{align*}
&\left\|w_\ell^i(t) - w_\ell^\star(t)\right\| \\
&\quad \leq A \rho^{t+1}\sum_{j=1}^M  \left\| w_\ell^j(0) \right\| + A \sum_{\tau = 0}^{t-1}  \sum_{j=1}^M \rho^{t-\tau} \varepsilon_{\tau+1}^j \Ind_{\left\{\tau \in T^j\right\}} \left\|H(\zbf_{\tau+1}^j,w^j(\tau))_\ell\right\| \\
& \quad \qquad (\mbox{by equation (\ref{eq:defSiDescent})})\\
& \quad \leq A \rho^{t+1}\sum_{j=1}^M  \left\| w_\ell^j(0) \right\| \\
& \qquad + A \sum_{\tau = 0}^{t-1} \sum_{j=1}^M \rho^{t-\tau} \varepsilon_{\tau+1}^j \Ind_{\tau \in T^j}\Ind_{\left\{\zbf^j_{\tau+1} \in W_\ell(w^j(\tau))\right\}} \left\|w_\ell^j(\tau) - \zbf^j_{\tau+1}\right\|.
\end{align*}
Therefore,
\begin{align*}
&\left\|w_\ell^i(t) - w_\ell^\star(t)\right\| \\
& \quad \leq A M \diam(\Gcal) \rho^{t+1} + A \diam(\Gcal) K_2 M \sum_{\tau = 0}^{t-1}\frac{1}{\tau \vee 1} \rho^{t-\tau}\\
 & \quad \qquad (\mbox{because $0 \in \Gcal$ and by Assumptions \ref{assum:strongdeacreasingSteps} and \ref{assum:witrajectories}})\\
 & \quad \leq A \diam(\Gcal) K_2  M\sum_{\tau = -1}^{t-1}\frac{1}{\tau \vee 1} \rho^{t-\tau}.\\
\end{align*}
Consequently,
\begin{align*}
&\left\|w^\star(t) - w^i(t) \right\| \\
& \quad = \sqrt{\sum_{\ell=1}^\kappa{\left\|w_\ell^i(t) - w_\ell^\star(t)\right\|^2}}\\
& \quad \leq \sqrt{\kappa}M\diam(\Gcal)A K_2\sum_{\tau=-1}^{t-1}\frac{1}{\tau \vee 1} \rho^{t- \tau}.
\end{align*}
This proves the desired result.
\proofend

Let us now introduce the following events: for any $\delta >0$ and $t\geq 0$,
\[A_\delta^t \triangleq \left\{w^\star(\tau) \in \Gcal_{\delta}^\kappa \;,\; \; t \geq \tau \geq 0\right\}.\]
Recall that the $\Gcal_{\delta}^\kappa$ is a compact subset of $\Gcal^\kappa$ defined by equality (\ref{eq:genCompact}). The next lemma establishes a detailed analysis of security regions for the parted components of the sequences $\left\{w^\star(t)\right\}_{t=0}^{\infty}$ and $\left\{w^j(t)\right\}_{t=0}^{\infty}$.
\begin{lemma}\label{lem:eventsproperties}
Let Assumptions \ref{assum:strongdeacreasingSteps} and \ref{assum:witrajectories} hold. Then,
\begin{enumerate}
\item there exists an integer $t^1_\delta \geq 1$ such that
 \[A_\delta^t \subset A_{\delta/2}^{t+1},\quad \mbox{$t \geq t^1_\delta$}.\]
 Moreover,
  \[w^\star(t) \in \Gcal_\delta^\kappa \Rightarrow [w^\star(t), w^\star(t+1)] \subset  \Gcal_{\delta/2}^\kappa, \quad \mbox{$ t \geq t^1_\delta$.}\]
\item There exists an integer $t^2_\delta \geq 1$ such that \\
\[w^\star(t) \in \Gcal_\delta^\kappa \Rightarrow [w^\star(t), w^i(t)] \subset  \Gcal_{\delta/2}^\kappa, \quad \mbox{$i \in\left\{1, \ldots, M \right\}$ and $t \geq t^2_\delta$.}\]
\end{enumerate}
\end{lemma}

{\bf Proof of Lemma \ref{lem:eventsproperties}.} {\bf Proof of statement 1.} The proof starts with the observation that under Assumption \ref{assum:witrajectories} we have $w^j(t) \in \Gcal^\kappa$, for all $i \in \left\{1, \ldots, M\right\}$ and $t \geq 0$. It follows that, for any $1 \leq \ell \leq \kappa$,
\begin{align*}
 \left\| H\left(\zbf_{t+1}^j,w^j(t)\right)_\ell \right\| & \leq \left\| \zbf_{t+1}^j - w_\ell^j(t) \right\| \\
& \leq \diam(\Gcal).
\end{align*}

Let us now provide an upper bound on the norm of the differences between two consecutive values of the agreement vector sequence. We may write, for all $t \geq 0$ and all $1 \leq \ell \leq M$,
\begin{align}
&\left\| w_\ell^\star(t + 1) - w_\ell^\star(t) \right\| \notag \\
& \quad = \left\| \sum_{j=1}^M{\phi^j(t)s_\ell^j(t)} \right\| \notag \\
& \quad \leq \sum_{j=1}^M{\phi^j(t)\left\|s_\ell^j(t) \right\|} \notag \\
& \quad \leq \sum_{j=1}^M{\varepsilon_{t+1}^j \Ind_{\left\{t \in T^j\right\}}\left\|H\left(\zbf_{t+1}^j,w^j(t) \right)_\ell \right\|} \notag \\
&  \quad \qquad (\mbox{by equation (\ref{eq:defSiDescent}) and statement 1. of Lemma \ref{lem:decompProp}}) \notag \\
& \quad \leq \frac{M \diam(\Gcal)K_2}{t \vee 1} \label{eq:ineqwstar}\\
&  \quad \qquad (\mbox{by Assumption \ref{assum:strongdeacreasingSteps}}). \notag
\end{align}

Take $t \geq \frac{4}{\delta}M \diam(\Gcal) K_2$ and $1 \leq k \neq \ell \leq M$. Let $\alpha$ be a real number in the interval $[0, 1]$.\\

If $w^\star(t) \in \Gcal^\kappa_\delta$ then
\begin{align*}
&\left\|(1 - \alpha ) w^\star_\ell(t) + \alpha w^\star_\ell(t+1) - (1-\alpha) w_k^\star(t) - \alpha w_k^\star(t+1) \right\| \\
& \quad = \left\| w^\star_\ell(t) -  w^\star_k(t) + \alpha\left( w^\star_\ell(t +1) - w^\star_\ell(t) \right)  + \alpha \left( w_k^\star(t) - w^\star_k(t+1)\right) \right\| \\
& \quad \geq \left\| w^\star_\ell(t) -  w^\star_k(t) \right\| - \left\| \alpha\left( w^\star_\ell(t +1) - w^\star_\ell(t)\right) + \alpha \left( w_k^\star(t) - w^\star_k(t+1)\right) \right\| \\
& \quad \geq \left\| w^\star_\ell(t) -  w^\star_k(t)\right\| - \alpha \left\|  w^\star_\ell(t +1) - w^\star_\ell(t)\right\|  - \alpha \left\|w_k^\star(t) - w^\star_k(t+1) \right\| \\
& \quad \geq \delta - 2\alpha \frac{\delta}{4}\\
& \quad \geq \delta/2.
\end{align*}
This proves that the whole segment $[w^\star(t), w^\star(t+1)]$ is contained in $\Gcal_{\delta/2}^\kappa$.\\

{\bf Proof of statement 2.} Take $t \geq 1$ and $1 \leq \ell \leq M$. If $w^\star(t) \in \Gcal^\kappa_\delta$ then
by Lemma \ref{lem:normMaj}, there exists $t_\delta^2$ such that
\[\left\| w_\ell^\star(t) - w_\ell^i (t) \right\| \leq \frac{\delta}{4}, \quad \mbox{$i \in \left\{1, \ldots, M\right\}$ and $t \geq t_{\delta}^2$.}\]
Let $k$ and $\ell$ two distinct integers between $1$ and $M$. For any $t \geq t_{\delta}^2$,
\begin{align*}
&\left\|\alpha w_k^i(t) +  (1 - \alpha)  w_k^\star(t) - \alpha w_\ell^i(t) - (1 - \alpha) w_\ell^\star(t) \right\| \\
& \quad = \left\| w_k^\star(t) -  w_\ell^\star(t) +\alpha(w_k^i(t) - w_k^\star(t)) + \alpha( w_\ell^\star(t) - w_\ell^i(t))\right\| \\
& \quad \geq \left\|  w_k^\star(t) -  w_\ell^\star(t) \right\| - \alpha \left\|w_k^i(t) - w_k^\star(t)\right\| -\alpha \left\|w_\ell^\star(t)  - w_\ell^i(t)\right\| \\
& \quad \geq \delta - 2 \alpha\frac{\delta}{4} \\
& \quad \geq \delta/2.
\end{align*}
This implies $[w^\star(t), w^i(t)] \subset  \Gcal_{\delta/2}^\kappa$, as desired.
\proofend

{\bf Proof of Proposition \ref{prop:epsilonstarspeed}.} By definition $\varepsilon_{t+1}^\star$ equals $\sum_{j=1}^M\Ind_{\left\{t \in T^j\right\}}\phi^j(t)\varepsilon_{t+1}^j$, for all $t \geq 0$ .\

On the one hand, since the real number $\phi^j(t)$ belongs to the interval $[\eta,1]$ (by Lemma \ref{lem:decompProp}) $\varepsilon_{t+1}^\star$ is bounded from above by $\frac{MK_2}{t \vee 1}$ using the right-hand side inequality of Assumption \ref{assum:strongdeacreasingSteps}.\\

On the other hand, $\varepsilon_{t+1}^\star$ is bounded from below by the nonnegative real number $\eta \frac{K_1}{t \vee 1 }$ using the left-hand side inequality of Assumption \ref{assum:strongdeacreasingSteps}. Note also that as Assumption \ref{assum:titimes} holds, this real number is a positive one. Therefore, it follows that
\[\varepsilon^\star_t \xrightarrow[t \rightarrow \infty]{} 0\]
and
\[\sum_{t = 0}^{\infty}{\varepsilon^\star_t } = \infty.\]
\proofend

{\bf Proof of Proposition \ref{prop:assumptionImplcations1}.} {\bf Consistency of $\sum_{t=0}^{\infty}{\Delta M_t^{(1)}}$.} Let $\delta$ be a positive real number and let $t \geq t_\delta^2$, where $t_\delta^2$ is given by Lemma \ref{eq:statement}. We may write
\begin{align*}
&\Ind_{A^t_\delta} \sum_{j=1}^M{\Ind_{\left\{t \in T^j\right\}} \phi^j(t)\varepsilon_{t+1}^j\left\|h\left(w^\star(t)\right)-h\left(w^j(t)\right) \right\|}\\
& \quad \leq \Ind_{\left\{[w^\star(t),w^j(t)] \subset \Gcal_{\delta/2 }^\kappa\right\}} \sum_{j=1}^M{
\phi^j(t)\varepsilon_{t+1}^j\left\|\nabla C\left(w^\star(t)\right)-\nabla C \left(w^j(t)\right) \right\|}\\
& \quad \qquad (\mbox{using statement 2. of Lemma \ref{lem:eventsproperties} and the fact that $\nabla C = h$ on $\Dcal_*^\kappa$}) \\
& \quad  \leq \Ind_{\left\{[w^\star(t),w^j(t)] \subset \Gcal_{\delta/2}^\kappa \right\}} P_{\delta/2} \sum_{j=1}^M{
 \varepsilon_{t+1}^j\left\|w^\star(t)-w^j(t) \right\|}\\
& \quad \qquad (\mbox{by Lemma \ref{lem:localLipschitz}})\\
  & \quad \leq \sqrt{\kappa}\diam(\Gcal)AK_2^2 P_{\delta/2}M^2 \frac{\theta_t}{t}\\
  & \quad \qquad (\mbox{by Lemma \ref{lem:normMaj}}).
\end{align*}
Thus, since $\sum_{t=0}^{\infty}{\frac{\theta_t}{t}} < \infty$, the series
\[ \sum_{t=0}^{\infty}{\Ind_{A^t_\delta} \sum_{j=1}^M{\Ind_{\left\{t \in T^j\right\}} \phi^j(t)\varepsilon_{t+1}^j\left\|h\left(w^\star(t)\right)-h\left(w^j(t)\right) \right\|}}\]
 is almost surely convergent. Under Assumption \ref{assum:wstartrajectories}, we have
 \[\mathbb{P}\left\{\bigcup_{\delta >0} \bigcap_{t \geq 0} A_\delta^t \right\} = 1.\]
It follows that the series $\sum_{t=0}^{\infty}{\Delta M_t^{(1)}}$ converges almost surely in $\rdkap$.\\

{\bf Consistency of $\sum_{t=0}^{\infty}{\Delta M_t^{(2)}}$.} The sequence of random variables $M^{(2)}_t$ defined,  for all $t\geq0$, by
\begin{align*}
M_t^{(2)} &\triangleq \sum_{\tau =0}^t \Delta M_\tau^{(2)}\\
&= \sum_{\tau=0}^t{\sum_{j=1}^M\Ind_{\left\{\tau \in T^j\right\}}\varepsilon_{\tau+1}^j \phi^j(\tau)\left(h\left(w^j(\tau)\right) - H\left(\zbf_{\tau+1}^j,w^j(\tau) \right)\right)}.
\end{align*}

is a vector valued martingale with respect to the filtration $\left\{\mathcal{F}_t\right\}_{t=0}^{\infty}$. It turns out that this martingale has square integrable increments. Precisely,

\[\sum_{t=0}^{\infty}{\esp\left\{\left\|M_{t+1}^{(2)} - M_t^{(2)} \right\|^2 \; {\Big \vert} \; \mathcal{F}_t \right\}}
 = \sum_{t=1}^{\infty}{\esp\left\{ \left\|\Delta M_t^{(2)}\right\|^2 \; {\Big \vert}  \; \mathcal{F}_t \right\}}
 < \infty.\]

Indeed, for all $j \in \left\{1 , \ldots, M\right\}$ and $t \geq 1$,
\begin{align*}
&\sum_{\tau = 1}^t{\esp\left\{\left\|\Ind_{\left\{\tau \in T^j\right\}}\varepsilon_{\tau+1}^j\left(h\left(w^j(\tau)\right) - H\left(\zbf_{\tau+1}^j(\tau),w^j(\tau) \right)\right) \right\|^2 \; {\big \vert} \; \mathcal{F}_\tau \right\}} \\
& \quad \leq \sum_{\tau = 1}^t{\left(\varepsilon_{\tau+1}^j\right)^2\esp\left\{\left\|h\left(w^j(\tau)\right) - H\left(\zbf_{\tau+1}^j(\tau),w^j(\tau)\right)\right\|^2 \; {\big \vert} \; \mathcal{F}_\tau \right\}} \\
& \quad \leq 2\sum_{\tau = 1}^t{\left(\varepsilon_{\tau+1}^j\right)^2\esp\left\{\left\|h\left(w^j(\tau)\right)\right\|^2 + \left\|H\left(\zbf_{\tau+1}^j(\tau),w^j(\tau)\right) \right\|^2 \; {\big \vert} \; \mathcal{F}_\tau \right\}} \\
& \quad \leq 4 \kappa \diam(\Gcal)^2 \sum_{\tau = 1}^t{\left(\varepsilon_{\tau+1}^j\right)^2}\\
& \quad \qquad (\mbox{using Assumption \ref{assum:witrajectories}})\\
& \quad \leq 4 \kappa \diam(\Gcal)^2 K_2^2 \sum_{\tau = 1}^t{\frac{1}{\tau^2}}.
\end{align*}

We conclude that the series $\sum_{t \geq 1} \Delta M^{(2)}_t$ is almost surely convergent.
\proofend

{\bf Proof of proposition \ref{prop:convSeries}.} Denote by $\langle x,y  \rangle$ the canonical inner product of two vectors $x,y \in \real^d$ and also, with a slight abuse of notation, the canonical inner product of two vectors $x,y \in \rdkap$. Let $\delta$ be a positive real number.\\
Take any $t\geq \max\left\{t_\delta^1,t_\delta^2 \right\}$, where $t_\delta^1$ and $t_\delta^2$ are defined as in Lemma \ref{lem:eventsproperties}. One has,
\begin{align*}
&\Ind_{A_\delta^{t+1}}C\left(w^\star(t+1)\right) \leq \Ind_{A_\delta^{t}}C\left(w^\star(t+1)\right). \\
& \quad \qquad (\mbox{by definition $A_\delta^{t+1} \subset A_\delta^{t}$}) \\
\end{align*}
Consequently,
\begin{align*}
&\Ind_{A_\delta^{t+1}}C\left(w^\star(t+1)\right) \\
& \quad \leq \Ind_{A_\delta^{t}}C\left(w^\star(t) \right) + \Ind_{A_\delta^{t}} \langle \nabla C(w^\star(t)) , w^\star(t+1) - w^\star(t) \rangle \\
& \qquad + \Ind_{\left\{[w^\star(t),w^\star(t+1)] \subset \Gcal^\kappa_{\delta/2}\right\}} \\
 & \quad \qquad  \times \left[\sup_{z \in \left[w^\star(t),w^\star(t+1) \right]}\left\{\left\|\nabla C(z) -\nabla C(w^\star(t))\right\|\right\}\left\| w^\star(t+1) - w^\star(t) \right\| \right] \\
& \quad \leq \Ind_{A_\delta^{t}}C\left(w^\star(t) \right) + \Ind_{A_\delta^{t}} \langle \nabla C(w^\star(t)) , w^\star(t+1) - w^\star(t) \rangle \\
& \qquad + P_{\delta/2}\left\| w^\star(t+1) - w^\star(t) \right\|^2 \\
& \quad \qquad (\mbox{using Lemma \ref{lem:localLipschitz}.}) \\
\end{align*}
The first inequality above holds since the bounded increment formula above is valid by statement 1 of Lemma \ref{lem:eventsproperties}. Let us now bound separately the right hand side members of the second inequality.\\

Firstly, the next inequality holds by inequality (\ref{eq:ineqwstar}) provided in the proof of Lemma \ref{lem:eventsproperties},
\[P_{\delta/2}\left\| w^\star(t+1) - w^\star(t) \right\|^2 \leq \kappa P_{\delta/2}\left(\frac{K_2 M \diam(\Gcal)}{t}\right)^2.\]
Secondly,
\begin{align*}
&\Ind_{A_\delta^{t}} \langle \nabla C(w^\star(t)) , w^\star(t+1) - w^\star(t) \rangle\\
&\quad = \Ind_{A_\delta^{t}}  \langle \nabla C(w^\star(t)) , \sum_{j=1}^M \phi^j(t)s^j(t) \rangle\\
&\quad \qquad (\mbox{by equation (\ref{eq:agreemRecursion})}) \\
& \quad = \Ind_{A_\delta^{t}}  \sum_{j=1}^M \langle \nabla C(w^j(t)) , \phi^j(t)s^j(t) \rangle \\
& \quad \qquad + \Ind_{A_\delta^{t}} \sum_{j=1}^M \langle\nabla C(w^\star(t)) - \nabla C(w^j(t)), \phi^j(t)s^j(t)\rangle.
\end{align*}
Thus,
\begin{align*}
&\Ind_{A_\delta^{t}} \langle \nabla C(w^\star(t)) , w^\star(t+1) - w^\star(t) \rangle\\
&\quad \leq \Ind_{A_\delta^{t}}\sum_{j=1}^M \langle \nabla C(w^j(t)) , \phi^j(t)s^j(t)\rangle \\
&\qquad + \Ind_{A_\delta^{t}} \sum_{j=1}^M \left| \langle \nabla C(w^\star(t)) -  \nabla C(w^j(t)) , \phi^j(t)s^j(t)\rangle \right|\\
&\quad \leq \Ind_{A_\delta^{t}}\sum_{j=1}^M \langle \nabla C(w^j(t)) , \phi^j(t)s^j(t)\rangle \\
& \qquad +  \sum_{j=1}^M \Ind_{A_\delta^{t}}\left\| \nabla C(w^\star(t)) -  \nabla C(w^j(t)) \right\| \left\| \phi^j(t)s^j(t) \right\|\\
& \quad \qquad (\mbox{using Cauchy-Schwarz inequality}).\\
\end{align*}
Therefore,
\begin{align*}
&\Ind_{A_\delta^{t}} \langle \nabla C(w^\star(t)) , w^\star(t+1) - w^\star(t) \rangle\\
&\quad \leq  \Ind_{A_\delta^{t}}\sum_{j=1}^M \langle \nabla C(w^j(t)) , \phi^j(t)s^j(t)\rangle \\
& \qquad +  \sum_{j=1}^M \Ind_{\left\{[w^\star(t),w^j(t)] \subset \Gcal_{\delta/2 }^\kappa\right\}} \left\| \nabla C(w^\star(t)) -  \nabla C(w^j(t)) \right\| \left\| \phi^j(t)s^j(t) \right\|\\
& \quad \qquad(\mbox{by statement 2 of Lemma \ref{lem:eventsproperties}})\\
& \quad \leq \Ind_{A_\delta^{t}} \sum_{j=1}^M \langle \nabla C(w^j(t)) , \phi^j(t)s^j(t)\rangle \\
& \qquad + P_{\delta/2} \sum_{j=1}^M  \left\| w^\star(t) - w^j(t) \right\| \left\| \phi^j(t)s^j(t) \right\|\\
& \qquad \quad (\mbox{using Lemma \ref{lem:localLipschitz}})\\
\end{align*}

\begin{align*}
&\Ind_{A_\delta^{t}} \langle \nabla C(w^\star(t)) , w^\star(t+1) - w^\star(t) \rangle\\
& \quad \leq \Ind_{A_\delta^{t}} \sum_{j=1}^M \langle \nabla C(w^j(t)) , \phi^j(t)s^j(t)\rangle \\
& \qquad + P_{\delta/2}A K_2^2 \kappa M^2 \diam(\Gcal)^2 \frac{\theta_t}{t}  \notag \\
& \qquad \quad (\mbox{using Lemma \ref{lem:normMaj} and the upper bound (\ref{eq:ineqwstar})}). \notag\\
\end{align*}

Finally,
\begin{align}
&\Ind_{A_\delta^{t+1}}C\left(w^\star(t+1)\right) \notag \\
& \quad \leq \Ind_{A_\delta^{t}}C\left(w^\star(t) \right) + \Ind_{A_\delta^{t}} \sum_{j=1}^M \langle \nabla C(w^j(t)) , \phi^j(t)s^j(t)\rangle \notag \\
& \qquad + P_{\delta/2}A K_2^2 \kappa M^2 \diam(\Gcal)^2 \frac{\theta_t}{t}  \notag \\
& \qquad + \kappa P_{\delta/2}\left(\frac{K_2 M \diam(\Gcal)}{t}\right)^2 \label{eq:technicalbound}.
\end{align}

Set
\[\Omega^1_\delta \triangleq  P_{\delta/2}A K_2^2 \kappa M^2 \diam(\Gcal)^2\]
and
\[\Omega^2_\delta \triangleq \kappa P_{\delta/2}\left(K_2 M \diam(\Gcal)\right)^2.\]

In the sequel, we shall need the following lemma.
\begin{lemma}\label{lem:supertechnicallemma}
For all $t \geq \max\left\{t_\delta^1,t_\delta^2 \right\}$, the quantity $W_t$ below is a nonnegative supermartingale with respect to the filtration $\left\{\mathcal{F}_t\right\}_{t=0}^{\infty}$:
\begin{align*}
W_t \triangleq \Ind_{A_\delta^t}C\left(w^\star(t)\right) & + \eta K_1\sum_{\tau = 0}^{t-1}{\Ind_{A_\delta^\tau}\frac{1}{\tau}\sum_{j=1}^M\Ind_{\left\{\tau \in T^j\right\}}\left\|\nabla C \left( w^j(\tau)\right) \right\|^2}\\
&+ \Omega^1_\delta \sum_{\tau = t}^{\infty}\frac{\theta(\tau)}{\tau} + \Omega^2_\delta \sum_{\tau = t}^{\infty}\frac{1}{\tau^2}, \quad t \geq 1.
\end{align*}
\end{lemma}
\newpage
{\bf Proof of Lemma \ref{lem:supertechnicallemma}.} Indeed, using the upper bound provided by equation (\ref{eq:technicalbound}),
\begin{align*}
&\esp \left\{\Ind_{A_\delta^{t+1}}C\left(w^\star(t +1 )\right) \; {\Big \vert} \; \mathcal{F}_t \right\}\\
& \quad \leq \Ind_{A_\delta^{t}}C\left(w^\star(t) \right) + \Ind_{A_\delta^{t}} \sum_{j=1}^M \esp\left\{\langle \nabla C(w^j(t)) , \phi^j(t)s^j(t)\rangle  \; {\big \vert} \; \mathcal{F}_t \right\} \\
& \qquad + \Omega^1_\delta \frac{1}{t}\theta_t + \Omega^2_\delta\frac{1}{t^2} \\
& \quad = \Ind_{A_\delta^{t}}C\left(w^\star(t) \right)\\
& \qquad + \Ind_{A_\delta^{t}} \sum_{j=1}^M \Big \langle \nabla C(w^j(t)) , \esp\left\{-\Ind_{\left\{t \in T^j\right\}}\phi^j(t)\varepsilon_{t+1}^j H(\zbf^j_{t+1},w^j(t))\rangle  \; {\Big \vert} \; \mathcal{F}_t \right\} \Big \rangle\\
& \qquad + \Omega^1_\delta \frac{\theta_t}{t} + \Omega^2_\delta\frac{1}{t^2} \\
& \quad = \Ind_{A_\delta^{t}}C\left(w^\star(t) \right)\\
& \qquad - \Ind_{A_\delta^{t}} \sum_{j=1}^M \Ind_{\left\{t \in T^j\right\}}\phi^j(t) \varepsilon_{t+1}^j \left\|\nabla C(w^j(t)) \right\|^2 + \Omega^1_\delta \frac{\theta_t}{t} + \Omega^2_\delta\frac{1}{t^2} \\
& \quad \leq \Ind_{A_\delta^{t}}C\left(w^\star(t) \right)\\
&  \qquad - \frac{\eta K_1}{t}\Ind_{A_\delta^{t}} \sum_{j=1}^M \Ind_{\left\{t \in T^j\right\}}\left\|\nabla C(w^j(t)) \right\|^2 + \Omega^1_\delta \frac{\theta_t}{t} + \Omega^2_\delta\frac{1}{t^2}.
\end{align*}
In the last inequality we used the fact that $\phi^j(t) \geq \eta$ (Lemma \ref{lem:decompProp}) and $\varepsilon_{t+1}^j \geq \frac{K_1}{t}$ (Assumption \ref{assum:strongdeacreasingSteps}).

It is straightforward to verify that, we have $W_t - \esp\{W_{t+1} \vert \mathcal{F}_t\} \geq0$ which prove the desired result.

\proofend

{\bf Proof of Proposition \ref{prop:convSeries} (continued).} Since $\left\{W_t\right\}_{t=1}^\infty$ is a nonnegative supermartingale (by Lemma \ref{lem:supertechnicallemma}), $W_t$ converges almost surely as $t \rightarrow \infty$ (see for instance Durrett \cite{DURR1}). Then, as $\sum_{\tau = t}^{\infty}\frac{\theta(\tau)}{\tau} \xrightarrow [t \rightarrow \infty]{} 0$ and $\sum_{\tau = t}^{\infty}\frac{1}{\tau^2} \xrightarrow[t \rightarrow \infty]{} 0$, we have
\begin{equation}\label{eq:convC}
\Ind_{A_\delta^t}C(w^\star(t)) \xrightarrow[t \rightarrow \infty]{} C_\infty, \quad \mbox{a.s.,}
\end{equation}
where $C_\infty \in [0, \infty)$
and, because the origin of the expression is increasing in $t$, the following series converges
\begin{equation}\label{eq:convSeries}
\sum_{\tau=0}^{\infty}{\Ind_{A_\delta^\tau}\frac{1}{\tau \vee 1}\sum_{j=1}^M\Ind_{\left\{\tau \in T^j\right\}}\left\|\nabla C \left( w^j(\tau)\right) \right\|^2}  < \infty, \quad \mbox{a.s.}
\end{equation}

{\bf Proof of statement 1.} Assumption \ref{assum:wstartrajectories} means that \[\mathbb{P}\left\{\bigcup_{\delta > 0} \bigcap_{t \geq 0} A_\delta^t\right\} =1.\]
Statement 1 follows easily from the convergence (\ref{eq:convC}).

{\bf Proof of statement 2.} The required convergence (\ref{eq:statement}) is proven as follows.
We have
\begin{align*}
&\sum_{\tau=0}^t{\varepsilon_{\tau+1}^\star \Ind_{A_\delta^\tau}\left\|\nabla C\left(w^\star(\tau)\right) \right\|^2}\\
& \quad \leq \sum_{\tau=0}^t{\sum_{j=1}^M\phi^j(\tau) \Ind_{\left\{\tau \in T^j\right\}}{\Ind_{A_\delta^\tau} \varepsilon_{\tau+1}^j\left\|\nabla C\left(w^\star(\tau)\right) \right\|^2}}\\
& \qquad \quad (\mbox{using equality (\ref{eq:defepsilonstar})})\\
& \quad \leq 2K_2\sum_{\tau=0}^t{\Ind_{A_\delta^\tau}\frac{1}{\tau \vee 1}\sum_{j=1}^M \Ind_{\left\{\tau \in T^j\right\}}{\left\|\nabla C\left(w^j(\tau)\right) \right\|^2}}\\
& \qquad \quad (\mbox{using Assumption \ref{assum:titimes}}) \\
& \qquad  + 2K_2\sum_{\tau=0}^t{\Ind_{\left\{[w^\star(\tau),w^j(\tau)] \subset \Gcal_{\delta/2 }^\kappa\right\}}\frac{1}{\tau \vee 1}\sum_{j=1}^M{\left\|\nabla C\left(w^j(\tau)\right) - \nabla C\left(w^\star(\tau)\right)  \right\|^2}}\\
& \qquad \quad (\mbox{using Assumption \ref{assum:titimes} and statement 2 of Lemma \ref{lem:eventsproperties}}.) \\
\end{align*}
Thus,
\begin{align*}
&\sum_{\tau=0}^t{\varepsilon_{\tau+1}^\star \Ind_{A_\delta^\tau}\left\|\nabla C\left(w^\star(\tau)\right) \right\|^2}\\
& \quad \leq 2K_2\sum_{\tau=0}^t{\Ind_{A_\delta^\tau}\frac{1}{\tau \vee 1}\sum_{j=1}^M \Ind_{\left\{\tau \in T^j\right\}}{\left\|\nabla C\left(w^j(\tau)\right) \right\|^2}}\\
& \qquad + 2K_2 P_{\delta/2}^2\sum_{\tau=0}^t{\Ind_{\left\{[w^\star(\tau),w^j(\tau)] \subset \Gcal_{\delta/2 }^\kappa\right\}}\frac{1}{\tau \vee 1}\sum_{j=1}^M{\left\|w^j(\tau) - w^\star(\tau)\right\|^2}} \\
& \qquad \quad (\mbox{by Lemma \ref{lem:localLipschitz}}).
\end{align*}

Thus,
\begin{align*}
& \sum_{\tau=0}^t{\varepsilon_{\tau+1}^{\star} \Ind_{A_{\delta}^\tau} \left\|\nabla C \left(w^\star(\tau)\right) \right\|^2}\\
& \quad \leq 2K_2\sum_{\tau=0}^t{\Ind_{A_\delta^\tau}\frac{1}{\tau \vee 1}\sum_{j=1}^M \Ind_{\left\{\tau \in T^j\right\}}{\left\|\nabla C\left(w^j(\tau)\right) \right\|^2}}\\
& \qquad + 2  P_{\delta/2}^2 K_2^3 \kappa M^3 A^2\diam(\Gcal)^2\sum_{\tau=1}^t{\frac{1}{\tau \vee 1 }\theta_\tau^2}\\
& \qquad \quad (\mbox{by Lemma \ref{lem:normMaj}}).
\end{align*}

Finally, using the convergence (\ref{eq:convSeries}), one has
\[ \sum_{\tau=0}^\infty{\varepsilon_{\tau+1}^\star \Ind_{A_\delta^\tau}\left\|\nabla C \left(w^\star(\tau)\right) \right\|^2} < \infty, \quad \mbox{a.s.},\]
and the conclusion follows from the fact that Assumption \ref{assum:wstartrajectories} implies
\[\prob\left\{\bigcup_{\delta > 0} \bigcap_{t \geq 0} A_\delta^t\right\} =1.\]
\proofend


{\bf Proof of Theorem \ref{thm:AsyncTheorem}.} The proof consists in verifying the assumptions of Theorem \ref{thm:AGLemma} with the function $\widehat{G}$ defined by equation (\ref{eq:Ghat}).\\

It has been outlined that Assumption \ref{assum:witrajectories} implies that $w^\star(t)$ lie in the compact set $\Gcal^\kappa$, almost surely, for all $t\geq 0$. Consequently, in the definition of $\widehat{G}(w^\star)$ the $\liminf$ symbol can be omitted. For all $\zbf \in \Gcal$ and all $t\geq0$, we have $\left\|H(\zbf,w^\star(t))\right\| \leq \sqrt{\kappa} \diam\left(\Gcal\right)$, almost surely, whereas $\left\{h(w^\star(t))\right\}_{t=0}^{\infty}$ satisfies
\[h(w^\star(t)) = \esp\left\{H\left(\zbf,w^\star(t)\right) \; \vert \; \mathcal{F}_t \right\}, \quad \mbox{$t \geq 0$, a.s.}\]
Thus, the sequences $\left\{w^\star(t)\right\}_{t=0}^{\infty}$ and $\left\{h(w^\star(t))\right\}_{t=0}^{\infty}$ are bounded almost surely.\\

Proposition \ref{prop:epsilonstarspeed}, respectively Proposition \ref{prop:assumptionImplcations1}, respectively Proposition \ref{prop:convSeries} show that the first assumption, respectively the third assumption, respectively the fourth assumption of Theorem \ref{thm:AGLemma} hold. This concludes the proof of the theorem.
\proofend

{\bf Justification of the statements (\ref{eq:thetaprop}).}
Recall that the definition of $\theta$ is provided in Lemma \ref{lem:normMaj}. Let us remark that it is sufficient to analyse the behavior in $t$ of the quantity $\sum_{\tau=1}^{t-1}{\rho^{t-\tau}/\tau}$.\\
Let $\varepsilon >0$ then for all $t \geq \lfloor 1/\varepsilon \rfloor +1$, we have
\begin{align*}
&\sum_{\tau=1}^{t-1}{\frac{\rho^{t-\tau}}{\tau}}\\
& \quad = \sum_{\tau=1}^{\lfloor 1/\varepsilon \rfloor }{\frac{\rho^{t-\tau}}{\tau}} + \sum_{\tau=\lfloor 1/\varepsilon \rfloor +1}^{t-1}{\frac{\rho^{t-\tau}}{\tau}} \\
& \quad \leq \sum_{\tau=1}^{\lfloor 1/\varepsilon \rfloor }{\rho^{t-\tau}} + \varepsilon \sum_{\tau=\lfloor 1/\varepsilon \rfloor +1}^{t-1}{\rho^{t-\tau}} \\
& \quad \leq \frac{\rho^{t-\lfloor 1/\varepsilon \rfloor}}{1-\rho} + \frac{\varepsilon}{1 - \rho} \\
&\quad \qquad (\mbox{using the fact that $\rho \in (0,1)$}).
\end{align*}
Consequently, for $t$ sufficiently large we have
\[\sum_{\tau=1}^{t-1}{\frac{\rho^{t-\tau}}{\tau}} \leq \frac{2\varepsilon}{1 -\rho}\]
which proves the first claim.\\
The second claim follows the same technique by letting ``$\varepsilon = 1/\sqrt{t}$''.\\
Thus, for $t \geq 1$ we have
\[\theta_t \leq \frac{\rho^{t - \lfloor \sqrt{t} \rfloor -1}}{1 - \rho} + \frac{1/\sqrt{t}}{1 - \rho}.\]
Finally, for $T\geq 1$, it holds
\[\sum_{t=1}^T{\sum_{\tau=1}^{t-1}{\frac{\rho^{t-\tau}}{\tau}}} \leq \frac{1}{1-\rho}\left(\sum_{t=1}^T{ \rho^{n-\lfloor \sqrt{n} \rfloor - 1} } + \sum_{t=1}^T{\frac{1}{n^{3/2}}} \right).\]
The two partial sums in the above parenthesis have finite limits which prove the second statement.

\end{document}